\documentclass[a4paper,14pt]{article}
\usepackage[paper=a4paper, dvips, top=3cm, left=2.7cm, right=2.7cm, foot=1.3cm, bottom=3.2cm]{geometry}
\usepackage{amsfonts,amssymb,amsmath,
amscd,amsthm,amsxtra,latexsym}
\usepackage{color,eufrak,setspace,
graphicx,xypic,fix-cm}

\newtheorem{thm}{Theorem}
\newtheorem{pro}[thm]{Proposition}
\newtheorem{lem}[thm]{Lemma}
\newtheorem{cor}[thm]{Corollary}
\theoremstyle{definition}
\newtheorem{dfn}[thm]{Definition}
\newtheorem{rem}[thm]{Remark}
\newtheorem{nta}[thm]{Notation}
\newtheorem{eg}[thm]{Example}

\title{Moduli of certain wild covers of curves}
\author{Jianru Zhang}
\date{}
\numberwithin{equation}{subsection}
\begin{document}
\maketitle
\begin{abstract}
A fine moduli space (see Section~\ref{secn&t} Definition~\ref{finemdli}) is constructed, for cyclic-by-$\mathsf{p}$ covers of an affine curve over an algebraically closed field $k$ of characteristic $\mathsf{p}>0$. An intersection (see Definition~\ref{M}) of finitely many fine moduli spaces for cyclic-by-$\mathsf{p}$ covers of affine curves gives a moduli space for $\mathsf{p}'$-by-$\mathsf{p}$ covers of an affine curve. A local moduli space is also constructed, for cyclic-by-$\mathsf{p}$ covers of $Spec(k((x)))$, which is the same as the global moduli space for cyclic-by-$\mathsf{p}$ covers of $\mathbb{P}^1-\{0\}$ tamely ramified over $\infty$ with the same Galois group. Then it is shown that a restriction morphism (see Lemma~\ref{res mor-2}) is finite with degrees on connected components $\textsf{p}$ powers: There are finitely many deleted points (see Figure 1) of an affine curve from its smooth completion. A cyclic-by-$\mathsf{p}$ cover of an affine curve gives a product of local covers with the same Galois group, of the punctured infinitesimal neighbourhoods of the deleted points. So there is a restriction morphism from the global moduli space to a product of local moduli spaces.
\end{abstract}

\textbf{Table of content}

1. Introduction

2. Notations and Terminology

3. Existence of moduli space for cyclic-by-$\mathsf{p}$ covers

4. Moduli for $\mathsf{p}'$-by-$\mathsf{p}$ covers

5. Local vs. global moduli

\section{Introduction} \label{intro}

The paper mainly generalizes the results in [H80] for $\mathsf{p}$-groups to cyclic-by-$\mathsf{p}$ groups defined in Section~\ref{secn&t} Definition~\ref{gpextn} a. See Section~\ref{secn&t} for notations and terminology below. Since [H80] is frequently cited, the statements of its main results are given in the Introduction.

In [H80], it is shown that (Theorem 1.2 [H80]) when $P$ is a finite $\mathsf{p}$-group, there exists a fine moduli space for pointed principal $P$-covers (see Section~\ref{secn&t} Definition~\ref{covdfn} b and c for the definition) of an affine curve over an algebraically closed field $k$ of characteristic $\mathsf{p}>0$, which is an ind affine space (see Definition~\ref{indschdfn} b). When $P'$ is a finite group whose order is prime to $\mathsf{p}$, there are only finitely many pointed principal $P'$-covers of an affine curve. The \textit{wild case}, where $\mathsf{p}$ divides the order of the Galois group of the cover, and the \textit{tame case}, where $\mathsf{p}$ does not divide the order of the Galois group of the cover, are very different. The fine moduli space for pointed principal $P$-covers of $\mathbb{P}^1-\{0\}$ gives a coarse moduli space for local pointed principal $P$-covers of $Spec(k((x)))$ (Proposition 2.1 [H80]). This is a special case of the next result, since the finite etale morphism there becomes an isomorphism here. Finally it is shown that a restriction morphism is finite etale (Proposition 2.7 [H80]), where the restriction morphism is described in the Abstract with the cyclic-by-$\mathsf{p}$ group there replaced by $P$ here. The result can be interpreted as a local-global principle: Given a pointed $P$-local cover at each of the deleted points of the affine curve from its smooth completion, there are only finitely many global pointed $P$-covers of the affine curve, whose restrictions at the deleted points are the ones given.
\begin{figure}
  \centering
  \includegraphics[width=4cm]{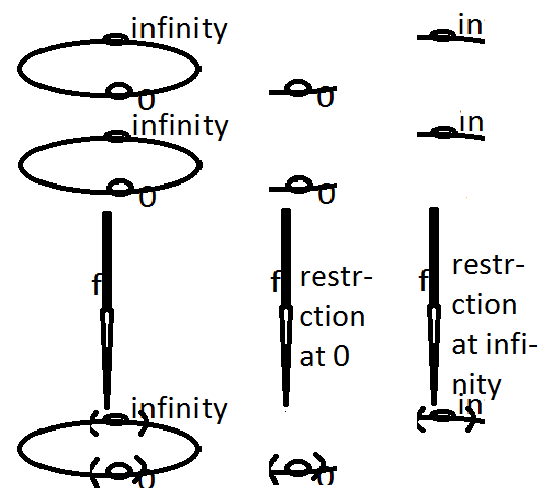}\\
  \caption{a trivial $\mathbb{Z}/2$-cover of $\mathbb{A}^1-\{0\}$}\label{pic1}
\end{figure}

In Figure 1, points 0 and $\infty$ are the deleted points of $\mathbb{A}^1-\{0\}$ from its smooth completion $\mathbb{P}^1$. The infinitesimal neighborhood of 0 is $Spec(k((x)))$ and the infinitesimal neighborhood of $\infty$ is $Spec(k((x^{-1})))$. $f$ gives a trivial $\mathbb{Z}/2$-cover of $\mathbb{A}^1-\{0\}$. The restriction of the global cover at 0 is a trivial $\mathbb{Z}/2$-cover of the infinitesimal neighborhood. Similarly for $\infty$.

The fine moduli space for pointed principal $P$-covers of an affine curve in [H80] is constructed in an inductive way with the base case for $P=\mathbb{Z}/\textsf{p}$.\\

Cyclic-by-$\mathsf{p}$ groups are the next simplest after $\mathsf{p}$-groups in the wild case. In the local situation, the Galois group of a connected Galois cover of $Spec(k((x)))$ is a cyclic-by-$\mathsf{p}$ group when $k$ is algebraically closed.

The fine moduli space for pointed principal cyclic-by-$\mathsf{p}$ covers of an affine curve is also constructed in an inductive way, using similar methods to those in the proof of Theorem 1.2 of [H80]. The fine moduli space is a disjoint union of finitely many ind affine spaces. Its relation to the the fine moduli space for pointed principal $P$-covers of an affine curve constructed in Theorem 1.2 [H80], is shown in Section~\ref{p' mdli} Lemma~\ref{(M)}.

The next simplest groups after cyclic-by-$\mathsf{p}$ groups are $\mathsf{p}'$-by-$\mathsf{p}$ groups defined in Section~\ref{secn&t} Definition~\ref{gpextn} a. A disjoint union of finitely many unions, of certain irreducible components in an intersection, of finitely many fine moduli spaces for cyclic-by-$\mathsf{p}$ covers of affine curves, gives a moduli space for $\mathsf{p}'$-by-$\mathsf{p}$ covers of an affine curve.

Two local-global principle results similar to those in [H80] described above are obtained, based on the construction of the fine moduli space for pointed principal cyclic-by-$\mathsf{p}$ covers of an affine curve, again using similar methods to those in [H80].\\

Here is the structure of the paper.

In Section~\ref{secn&t}, notations and terminology are given, which are used throughout Sections 3, 4, and 5 without explanation again. In Section~\ref{glb mdli}, a fine moduli space for pointed principal $G$-covers of an affine curve (Theorem~\ref{thm1}), where $G$ is a cyclic-by-$\mathsf{p}$ group, is constructed. In Section~\ref{p' mdli}, it is shown that a disjoint union of finitely many unions, of certain irreducible components in an intersection, of finitely many fine moduli spaces for cyclic-by-$\mathsf{p}$ covers of some affine curves, gives a moduli space for $\mathsf{p}'$-by-$\mathsf{p}$ covers of an affine curve (Corollary~\ref{uniirr}). In Section~\ref{l vs g}, a global fine moduli space is constructed (Proposition~\ref{ram mdli-3}) for cyclic-by-$\mathsf{p}$ covers of an affine curve at most tamely ramified over finitely many closed points, as well as a parameter space for local cyclic-by-$\mathsf{p}$ covers of $Spec(k((x)))$ (Proposition~\ref{local para space-3}). Then it is shown that a restriction morphism is finite with degrees on connected components $\textsf{p}$ powers, which is from the global moduli space to a product of the local parameter spaces (Proposition~\ref{fin etl-1}).

Leitfaden:
\[\xymatrix{
  Section~\ref{secn&t}\ar@{.>}[d]|-{}  \\
  Section~\ref{glb mdli} \ar@/^/[drr]^{}
    \ar@/_/[dr]|-{}                   \\
   & Section~\ref{p' mdli}  & Section~\ref{l vs g}
   .}
\]

Similar work can be found in [K86] and [P02]. In [K86], Main Theorem 1.4.1 is essentially the version over a general field of characteristic $\mathsf{p}>0$ of Proposition~\ref{local-global moduli}. In [P02], a configuration space $C(I,j)$ is constructed in Section 2.2, which is for $I=\mathbb{Z}/\mathsf{p}\rtimes\mathbb{Z}/n$-covers of $Spec(k[[u^{-1}]])$ with jump $j$. This is related to the parameter space given in Proposition~\ref{local para space-3}.\\

For the results in [H80] below, $Gr$ is a finite group, $P$ a finite $\textsf{p}$-group, $k$ an algebraically closed field, $(U_0,u_0)$ geometrically pointed $Spec(k((x)))$ and $(U,u_g)$ a geometrically pointed affine curve, as defined in Section~\ref{general setting}.

Let $(S,s_0)$ be a pointed (see Definition~\ref{schs} b and Definition~\ref{pted} a) connected affine $k$-scheme. In Definition~\ref{pted} b, when two pointed $S$-parameterized $Gr$-covers of $U$ are equivalent is defined; two such covers are equivalent if they agree pulled back to a finite etale cover $T$ of $S$. Since a pointed $S$-parameterized $Gr$-cover of $U$ corresponds to a homomorphism $\tilde{\varphi}:\pi_1(S\times U,(s_0,u_g))\rightarrow Gr$, the definition of an equivalence class of $\tilde{\varphi}$ is induced in Remark~\ref{fml} a. Similarly for the local case; in Definition~\ref{w}, the $w$-equivalence class of a pointed $S$-parameterized $Gr$-cover of $Spec(k((x)))$ is defined, which induces the definition of a $w$-equivalence class of a homomorphism $\tilde{\varphi}:\pi_1(S\times Spec(k((x))),(s_0,u_0))\rightarrow Gr$.

With these terminology, the following definition can be given.
\begin{dfn}\label{H80F}
a. Define $F_{U,P}$ as the functor: $\mathcal{S}_1\rightarrow$(Sets); $(S,s_0)\mapsto \{[\tilde{\varphi}]$, where $\tilde{\varphi}:\pi_1(S\times U,(s_0,u_g))\rightarrow P$ is a group homomorphism and $[\tilde{\varphi}]$ is the equivalence class (see above or Remark~\ref{fml} a) of $\tilde{\varphi}$.\}

b. Define $F_{U_0,P}^{w}$ as the functor: $\mathcal{S}_1\rightarrow$(Sets); $(S,s_0)\mapsto \{[\tilde{\varphi}]^w$, where $\tilde{\varphi}:\pi_1(S\times U_0,(s_0,u_0))\rightarrow P$ is a group homomorphism and $[\tilde{\varphi}]^w$ is the w-equivalence class (see above or Definition~\ref{w}) of $\tilde{\varphi}$.\}

\end{dfn}
See Section~\ref{fine moduli space} for the definition of a fine moduli space. Theorem~\ref{thm1.2} means that $M_{U,P}$ represents the moduli functor $F_{U,P}$. The ``direct limit of affine spaces" in Theorem~\ref{thm1.2} is called an ind affine space in the paper (See Section~\ref{ind schemes} Definition~\ref{indschdfn} b).

\begin{thm}\label{thm1.2}
(Theorem 1.2, [H80]) There is a fine moduli space $M_{U,P}$ (denoted by $M_G$ there with $G=P$) for pointed families of principal $P$-covers of $U$, namely a direct limit of affine spaces $\mathbb{A}_k^N$.
\end{thm}

The local case, moduli problem for $P$-covers of $Spec(k((x)))$, is simpler than the global case. A parallel construction to the one in the global case gives a coarse moduli space of pointed $P$-covers of $Spec(k((x)))$. Proposition~\ref{pro2.1} below means that $M_{\mathbb{P}^1-\{0\},P}$ represents the moduli functor $F^w_{U_0,P}$.
\begin{pro}\label{pro2.1}
(Proposition 2.1, [H80]) The fine moduli space $M_{\mathbb{P}^1-\{0\},P}$ for pointed principal $P$-covers of $\mathbb{P}^1-\{0\}$ is also a coarse moduli space for pointed principal $P$-covers of $Spec(k((x)))$, compatibly with the inclusion $Spec(k((x)))\subseteq \mathbb{P}^1-\{0\}$.
\end{pro}

\begin{pro}\label{pro2.7}
(Proposition 2.7, [H80]) Let $M_{U,P}\rightarrow\Pi_i M^l_{U_{0_i},P}$ be the restriction morphism described in the Abstract. It is an etale cover. Its degree is a power of $\textsf{p}$, and is equal to the number of pointed principal $P$-covers of the completion $\bar U$.
\end{pro}

\begin{thm}\label{thm1.12}
(Theorem 1.12, [H80]) Let $M_{U,P}$ be the fine moduli space for pointed principal $P$-covers of $(U,u_g)$. There is a natural action of $Aut(P)$ on $M_{U,P}$, and a dense open subset $M^0_{U,P}$ of $M_{U,P}$ parameterizing connected principal covers, such that $\bar M^0_{U,P}:=M^0_{U,P}/Aut(P)$ is a fine moduli space for pointed families of Galois covers of $(U,u_g)$ with group $P$.
\end{thm}

\textbf{Acknowledgement} I thank D. Harbater for helpful discussions during the preparation of the paper.

\section{Notations and Terminology}\label{secn&t}

Terms and symbols are defined here that will be used in Sections 3, 4, and 5 without being explained again.
\subsection{General settings}\label{general setting}
\begin{dfn}\label{gpextn}
a. Groups of the form $P\rtimes_{\rho} \mathbb{Z}/n$ are called \textit{cyclic-by-$\mathsf{p}$ groups}, where $\mathsf{p}$ is a prime number, $P$ a finite $\mathsf{p}$-group and $\rho:\mathbb{Z}/n \rightarrow Aut(P)$ an action of $\mathbb{Z}/n$ on $P$ with $n$ and $\mathsf{p}$ coprime. Groups of the form $P\rtimes_{\rho'} P'$ are called \textit{$\mathsf{p}'$-by-$\mathsf{p}$ groups}, where $P'$ is a finite group whose order is prime to $\mathsf{p}$ and $\rho': P'\rightarrow Aut(P)$ an action of $P'$ on $P$.

b. Let $n_t$ be a factor of $n$, $x_t=n/n_t$, $\iota_{n_t}$ the embedding of $\mathbb{Z}/n_t$ into $\mathbb{Z}/n$ sending $\bar 1\in \mathbb{Z}/n_t$ to $\bar x_t\in \mathbb{Z}/n$ and $\rho_{n_t}=\rho\circ \iota_{n_t}$.

c. $Gr$ always represents an arbitrary finite group and $G$ represents $P\rtimes_{\rho} \mathbb{Z}/n$.
\end{dfn}

\begin{dfn}\label{schs}
a. Let $k$ be an algebraically closed field of characteristic $\mathsf{p}>0$ and fix a primitive $n$-th root of unity $\zeta_n$ in $k$. Write $U_0=Spec(k((x)))$ and $\bar U_0=Spec(k[[x]])$. Denote the fiber product $S\times _kX$ by $S\times X$, where $S$ and $X$ are $k$-schemes.

b. \textit{Pointed} means geometrically pointed unless otherwise stated. A geometric point of a scheme $X$ is a morphism from $Spec(\Omega)$ to $X$ with $\Omega$ an algebraically closed field. \textit{Curve} means a connected smooth integral affine 1 dimensional scheme of finite type over $k$.

c. Denote by $\mathcal{S}$ (resp. $\mathcal{S}_1$) the full subcategory of the category (Pointed $k$-schemes) of all pointed $k$-schemes, whose objects are connected affine pointed finite type $k$-schemes (resp. connected affine pointed $k$-schemes). Denote by $\mathcal{S}'$ (resp. $\mathcal{S}_1'$) the non pointed version of $\mathcal{S}$ (resp. $\mathcal{S}_1$).

d. $(U,u_g)$ always represents a pointed curve.
\end{dfn}

\begin{rem}\label{lift}
The word ``lift" has two meanings in the paper. The first meaning is to extend a group homomorphism whose domain is a fundamental group, to a group homomorphism with domain a bigger fundamental group. This meaning is given around diagram (\ref{glb mdli}.1) in the definition of $\rho$-liftable. The second meaning is to lift a morphism $ \bar{\phi}$ mapping to a quotient group $\bar P$, to some morphism $\phi$ mapping to the original group $P$:
\[\xymatrix{
                &         P \ar[d]^{}     \\
  \pi_1(U,u_g) \ar[ur]^{\phi} \ar[r]_{\bar{\phi}} & \bar P             .}
\]
\end{rem}
\subsection{Covers}
\begin{dfn}\label{covdfn}
a. Let $X$ be a connected scheme and ($Fet_X$) the category of finite etale covers of $X$. For every point $x$ of $X$ (recall from Definition~\ref{schs} b that a point on $X$ always means a geometric point), the \textit{fiber functor} $Fib_x$: ($Fet_X$)$\rightarrow$(Sets) sends a finite etale cover $Y\xrightarrow{pr}X$ to $Y_x$, the geometric fiber (pullback of $Y$ to $x$) of $Y$ at $x$. The fundamental group $\pi_1(X,x)$ is defined as the automorphism group of the fiber functor $Fib_x$. Then $Y_x$ is a left $\pi_1(X,x)$-set.

b. A \textit{principal $Gr$-cover} of a connected scheme $X$ not necessarily over $k$, is a finite etale cover $Z\rightarrow X$ together with an embedding of $Gr$ in the group $Aut(Z/X)$, such that $Gr$ acts simply transitively (left group action) on every geometric fiber of $Z\rightarrow X$. A \textit{$Gr$-cover} means a principal $Gr$-cover.

c. With the notations of b, $Z$ is pointed over $x_0$ means $Z$ is pointed at some point $z_0$ that maps to $x_0$ under $Z\rightarrow X$. Two pointed $Gr$-covers of $(X,x_0)$ are isomorphic if there is an isomorphism between them:
\[\xymatrix{
  (Z,z_0) \ar[rr]^{f\simeq} \ar[dr]_{}
                & &  (Z',z'_0)    \ar[dl]^{}    \\
                & (X,x_0) ,
                }\]
such that the triangle diagram commutes and the diagram
\[\xymatrix{
  Z \ar[d]_{g} \ar[r]^{f} & Z' \ar[d]^{g} \\
  Z \ar[r]^{f} & Z'   }\]
commutes for each $g\in Gr$.
\end{dfn}

\begin{rem}\label{abel}
a. There is a natural bijection between the set of isomorphism
classes of $Gr$-covers of $X$ pointed over a fixed base point
$x_0$, and the set of homomorphisms from $\pi_1(X,x_0)$ to
$Gr$, hence below a pointed $Gr$-cover is often identified with
the homomorphism corresponding to it.

b. If $Gr$ is abelian, the
set of homomorphisms from $\pi_1(X,x_0)$ to
$Gr$ is a group. It may be identified with the etale cohomology group $H^1(X,Gr)$ (SGA 1, XI 5, or Example 11.3 in [M13]), or in terms of group cohomology with $H^1(\pi_1(X,x_0),Gr)$.

c. Let $Gr$ be an abelian group and $W\rightarrow U$ be a $Gr$-cover of $U$. Then $W$ pointed at a point $w_g$ over $u_g$ is isomorphic to, as pointed $Gr$-covers of $(U,u_g)$, $W$ pointed at any other point $w'_g$ over $u_g$.
\end{rem}

\begin{dfn}\label{pted}
a. A point $(s_0,v_g)$ on a fiber product $S\times V$ of $k$-schemes means a commutative diagram by the universal property of a fiber product:
\[
 \xymatrix{
  Spec(\Omega) \ar[d]_{s_0} \ar[r]^{v_g} & V\ar[d]^{} \\
  S \ar[r]^{} & Spec(k) ,}\eqno{(\ref{secn&t}.1)}
 \]
where $\Omega$ is some algebraically closed field.

b. A \textit{pointed family of $Gr$-covers} of a pointed connected $k$-scheme $X$, parametrized by a pointed connected affine $k$-scheme $S$, means an equivalence class of pointed $Gr$-covers of $S\times X$, two being equivalent if they become isomorphic after being pulled back by some finite etale cover $(T,t_0)\rightarrow (S,s_0)$.
\end{dfn}

\begin{rem}\label{fml}
a. Two elements $\widetilde{\phi}$ and $\widetilde{\phi}'$ in $Hom(\pi_1(S\times U, (s_0,u_g)),Gr)$ are \textit{equivalent} if their corresponding pointed $Gr$-covers of $(S\times U, (s_0,u_g))$ are equivalent. Denote the equivalence class of $\widetilde{\phi}$ by $[\widetilde{\phi}]$.

b. Using equivalence classes (see Definition~\ref{F1}, Definition~\ref{F1pair}, Definition~\ref{F2}, Definition~\ref{F3}), rather than isomorphism classes, a fine moduli space can be constructed. The definition of equivalence, using finite etale covers, arises naturally in the proof of Theorem ~\ref{thm3}: Assume for instance $e_{\rho}=1$, then the last paragraph in the proof gives $F(S,s_0)=\frac{H^1(S\times U,\mathbb{F}_q)}{H^1(S,\mathbb{F}_q)}$. The equality holds because the  definition of equivalence uses finite etale covers (see c). With the equality it can be shown that $F$ is represented by an ind affine space.

c. If $Gr$ is abelian, then the set of such pointed families may be identified with $H^1(S\times X,Gr)/H^1(S,Gr)$, where $H^1(S\times X,Gr)$ and $H^1(S,Gr)$ are standard etale cohomology groups.
\end{rem}
\begin{dfn}\label{chemindfn}
Suppose $X$ is a connected scheme and $x,x'$ two geometric points on $X$. A \textit{chemin} $x'\rightarrow x$ means an isomorphism from the fiber functor $Fib_x$ to the fiber functor $Fib_{x'}$ ([S09], Remark 5.5.3, p171). Since the fundamental group $\pi_1(X,x)$ is defined as the automorphism group of the fiber functor $Fib_x$, a chemin $x'\rightarrow x:Fib_x\xrightarrow{i}Fib_{x'}$ induces an isomorphism $\pi_1(X,x)\xrightarrow{\simeq}\pi_1(X,x'):
\alpha\mapsto i\alpha i^{-1}$.
\end{dfn}

\subsection{Ind schemes}\label{ind schemes}
\begin{dfn}\label{indschdfn}
a. An \textit{ind scheme} means, in the paper, a direct system of $k$-schemes $\{X_i\}$ indexed by natural numbers with transition $k$-morphisms \{$X_i\xrightarrow{x_i} X_{i+1}$\}.

b. An ind scheme is an \textit{ind affine space}, if every $X_i$ is an affine space $\mathbb{A}_k^{n_i}$.
\end{dfn}
\begin{rem}\label{compo}
Every moduli space $M$ in the paper is a disjoint union of finitely many ind affine spaces, then each ind affine space is called \textit{a connected component} of $M$. An ind affine space $M$ can be viewed as a functor: $\mathcal{S}_1\rightarrow$(Sets); $(S,s_0)\mapsto\;Hom(S,M)$. The disjoint union of finitely many ind affine spaces $\{M_i\}$ is the functor $\amalg M_i$: $\mathcal{S}_1\rightarrow$(Sets); $(S,s_0)\mapsto\;\amalg_i Hom(S,M_i)$.
\end{rem}
\begin{dfn}\label{indschmorphdfn}
a. A \textit{pre-morphism} from a $k$-scheme $X$ to an ind scheme $\{X_m\}$ is the equivalence class of a $k$-morphism between schemes $g_{m_0}:X\rightarrow X_{m_0}$ for some $m_0$, where two morphisms $g_{m_0}$ and $g_{m_1}$ are equivalent if for some $m_2\geq m_0,m_1$ the two composition morphisms $X\xrightarrow{g_{m_0}}X_{m_0}\rightarrow X_{m_2}$ and $X\xrightarrow{g_{m_1}}X_{m_1}
\rightarrow X_{m_2}$ are the same.

b. A \textit{pre-morphism} from an ind scheme $\{X_m\}$ with transition morphisms $\{x_m\}$ to another ind scheme $\{Y_m\}$ with transition morphisms $\{y_m\}$, is the equivalence class of a system of compatible $k$-morphisms $\{f_m|m\geq m_0\}$ between schemes with $f_m:X_m\rightarrow Y_{N_m}$. The system $\{f_m|m\geq m_0\}$ is compatible means that for every $m\geq m_0$, there exists an $n_m$ such that the following diagram is commutative:
\[\xymatrix{
  X_m \ar[d]_{x_m} \ar[r]^{f_m} & Y_{N_m}  \ar[dr]^{y_{N_mn_m}}\\
  X_{m+1} \ar[r]^{f_{m+1}} & Y_{N_{m+1}}\ar[r]^{y_{N_{m+1}n_m}}
  &Y_{n_m},}
\]
where $y_{N_mn_m}$ is the transition morphism from $Y_{N_m}$ to $Y_{n_m}$ and similarly for $y_{N_{m+1}n_m}$. Two compatible systems $\{f_m|m\geq m_0\}$ and $\{g_m|m\geq m_1\}$ are \textit{equivalent}, if there exists an $m_2\geq m_0,m_1$ such that for every $m\geq m_2$ the two morphisms $f_m$ and $g_m$ are equivalent in the sense of a. Every pre-morphism between ind schemes in the paper, in Lemma~\ref{(M)}, Lemma~\ref{cmap}, Lemma~\ref{res mor-1} and Lemma~\ref{res mor-2}, can be given by a compatible system $\{f_m\}$ of the special form: $X_m\xrightarrow{f_m} Y_m$ and the following diagram commutes
\[\xymatrix{
  X_m \ar[d]_{x_m} \ar[r]^{f_m} & Y_{m}  \ar[d]^{y_m}\\
  X_{m+1} \ar[r]^{f_{m+1}} & Y_{m+1}.}(\ref{indschmorphdfn}.1)
\]

c. In either a or b, a presheaf can be gotten. In b, the presheaf $Pre$ is from the site of (ind schemes)$\times$(ind schemes) with etale topology to (sets); $(\{X_m\},\{Y_m\})\mapsto PreMorph(\{X_m\},\{Y_m\})$. Let $sPre$ be the sheafification of $Pre$. A \textit{morphism} between $\{X_m\}$ and $\{Y_m\}$ is an element in $sPre(\{X_m\},\{Y_m\})$. Similarly for a. Since the construction is canonical, it suffices to check assertions on pre-morphisms.

d. In the cases of Lemma~\ref{(M)}, Lemma~\ref{cmap}, Lemma~\ref{res mor-1} and Lemma~\ref{res mor-2}, the morphism given by $\{f_m\}$ in diagram (\ref{indschmorphdfn}.1), is \textit{surjective} (resp. \textit{finite, finite etale}), if there exists some natural number $m_0$ such that for every $m\geq m_0$ the $k$-morphism $f_m$ is surjective (resp. finite, finite etale).
\end{dfn}
\subsection{Fine moduli space}\label{fine moduli space}
\begin{dfn}\label{finemdli}
A \textit{fine moduli space} $M$ for a contravariant functor $F$ from the category $\mathcal{S}_1$ to the category (Sets), is an ind scheme such that $F$ is isomorphic to the functor $Hom(\bullet,M)$: $\mathcal{S}_1\rightarrow$(Sets); $(S,s_0)\mapsto$ \{$k$-morphisms from $S$ to $M$\}.
\end{dfn}

Below is a list of moduli functors in the paper, with their rough meanings and places where they are defined.

\textbf{List of moduli functors}

$F_{U,P}$; the functor for pointed $P$-covers of $(U,u_g)$ using equivalence classes; Definition~\ref{H80F} a

$F_{U_0,P}^{w}$; the functor for pointed $P$-covers of $(Spec(k((x))),u_0)$ using w-equivalence classes; Definition~\ref{H80F} b

$F_{V,H}^{\rho}$; the functor for pointed $\rho$-liftable $H$-covers of $(V,v_g)$ using equivalence classes; Definition~\ref{F1} and Definition~\ref{newlftb}

$F_{V,H}^{\rho\bullet}$; the functor for $\rho$-liftable pairs (of pointed $H$-covers) of $(V,v_g)$ using equivalence classes; Definition~\ref{F1pair}

$F_{V,P}^{\rho\bullet}$; the functor for $\rho$-liftable pairs (of pointed $P$-covers) of $(V,v_g)$ using equivalence classes; Definition~\ref{F2}

$F_{U,G}$; the functor for pointed $G$-covers of $(U,u_g)$ using equivalence classes; Definition~\ref{F3}

$F_{U,G}^T$; the functor for pointed $G$-covers of $(U-T,u_g)$, at most tamely ramified over $T$ consisting of finitely many closed points on $U$, using equivalence classes; Definition~\ref{F1T}

$F_{V_l,P}^{\rho_{n_l}\bullet/T}$; the functor for $\rho_{n_l}$-liftable pairs (of pointed $P$-covers) of $(V_l,v_l)$ using equivalence classes, with $V_l$ a cover of $U$ at most ramified over $T$ consisting of finitely many closed points on $U$; Definition~\ref{F2T}

$F_{V_0,P}^{w\rho\bullet}$; the functor for $\rho$-liftable pairs (of pointed $P$-covers) of $(V_0,v_0)$ using w-equivalence classes; Definition~\ref{w}\\

Comments about several subtle concepts are collected below for reference convenience.

Comments concerning ind schemes include Remark~\ref{compo}.

Comments concerning universal families include Remark~\ref{uni fml}, Definition~\ref{unifmlpair}, Remark~\ref{unifml2}, Remark~\ref{n}, and Definition~\ref{attach}.

Comments concerning fine moduli spaces include Remark~\ref{fml}, Remark~\ref{clpt}, Remark~\ref{compo}, Remark~\ref{n}, Remark~\ref{dcps}, Remark~\ref{clsub}, and  Definition~\ref{attach}.
\subsection{Table of symbols}\label{table of symbols}
Below is a table of symbols, which are used in Sections 3, 4 and 5 without explanation again after their definitions. It gives meanings of symbols and places where they are defined. ``Beginning" of a section means beginning of the rest of the section below the introductory part.

\begin{center}
\textbf{Table of symbols}\\
\end{center}

$c$; a fixed element in $\pi_1(U,u_g)$ that maps to $\bar{1}$
under $\theta$; Section~\ref{glb mdli}, beginning

$c_i$: similar to $c$

$c'_i$; a fixed element in $\pi_1(V_i',\overline{v'_{gi}})$ that maps to $p_i'$
under $\theta_i'$; Section~\ref{p' mdli}, beginning

$H$; an elementary abelian group of order a $\textsf{p}$-power;
 Lemma~\ref{reg}

$\{pr_i:(V_i,v_i)\rightarrow(U,u_g)\}$; the set of all connected
pointed $\mathbb{Z}/n_i$-covers of $(U,u_g)$ with $n_i$
running over factors of $n$; Section~\ref{glb mdli}, beginning

$T$; a finite set of closed points on $U$ not including $u_g$; Section~\ref{l vs g}, beginning

$U^0$; $U-T$; Section~\ref{l vs g}, beginning

$(U_0,u_0)$; pointed $Spec(k((x)))$; Section~\ref{l vs g}, Notation~\ref{V0t}

$(V,v_g)$; a fixed connected pointed $\mathbb{Z}/n$-cover of
$(U,u_g)$; Section~\ref{glb mdli}, beginning

$(V',v'_g)$; a fixed connected pointed $P'$-cover of
$(U,u_g)$; Section~\ref{p' mdli}, beginning

$(V'_{i'},\overline{v'_{gi}})$; quotient of $(V',v'_g)$ by $\langle p'_i\rangle$; Section~\ref{p' mdli}, beginning

$\{(V_l^0,v_l)\}$; the set of all connected
pointed $\mathbb{Z}/n_l$-covers of $(U^0,u_g)$ with $n_l$
running over factors of $n$; Section~\ref{l vs g}, beginning

$V_l$; extension of $V_l^0$, by putting back in the closed points over $T\subset U$ to $V_l^0$, which are originally missing from $V_l^0$'s smooth completion; Section~\ref{l vs g}, beginning

$(V_{0t},v_{0t})$; a connected pointed $\mathbb{Z}/n_t$-cover of $(U_0,u_0)$ with $n_t$ a factor of $n$; Section~\ref{l vs g}, Notation~\ref{V0t}

$\rho$; an action of $\mathbb{Z}/n$ on $P$; Section~\ref{general setting} Definition~\ref{gpextn} a

$\rho'$; an action of $P'$ on $P$; same as $\rho$ above

$[\widetilde{\phi}]$; the equivalence class of $\widetilde{\phi}$;
Remark~\ref{fml}

$\rho'_i$; an action of $\langle p_i'\rangle$ on $P$ given by
restriction of $\rho'$; Section~\ref{p' mdli} below Remark~\ref{M's}

$\theta$; the group homomorphism $\pi_1(U,u_g)\rightarrow
\mathbb{Z}/n$ corresponding to $(V,v_g)\rightarrow(U,u_g)$;
Section~\ref{glb mdli}, beginning

$\theta_i$: similar to $\theta$

$\theta'$; the group homomorphism $\pi_1(U,u_g)\rightarrow
P'$ corresponding to $(V',v'_g)\rightarrow(U,u_g)$;
Section~\ref{p' mdli}, beginning

$\theta_i'$; the group homomorphism $\pi_1(V'_i,\overline{v'_{gi}})
\rightarrow
\langle p_i'\rangle$ corresponding to $(V',v'_g)\rightarrow(V'_i,
\overline{v'_{gi}})$;
Section~\ref{p' mdli}, beginning

\section{Existence of moduli space for cyclic-by-$\mathsf{p}$ covers}\label{glb mdli}
In Section~\ref{glb mdli}, a fine moduli space that represents the functor $F_{U,G}$ defined above Theorem~\ref{thm1}, for pointed $G$-covers of the pointed affine curve $(U,u_g)$, where $G$ is a cyclic-by-$\mathsf{p}$ group, is constructed. The construction is done in 3 steps: Theorem~\ref{thm3}$\Rightarrow$Theorem~\ref{thm2}$\Rightarrow$
Theorem~\ref{thm1}. Theorem~\ref{thm3} is the base case of an induction, and Theorem~\ref{thm2} is the inductive step. Theorem~\ref{thm1} collects building blocks given in Theorem~\ref{thm2} to build the target fine moduli space.

As always, we follow notations and terminology defined in Section~\ref{secn&t}. For example $G$ represents a cyclic-by-$\mathsf{p}$ group.\\

Here are necessary settings for Theorem~\ref{thm3}.

Since $(n,\mathsf{p})=1$, for every factor $n'$ of $n$, there are only finitely many connected pointed $\mathbb{Z}/n'$-covers of $(U,u_g)$, up to isomorphism. See also Remark~\ref{abel}.

Denote these covers by $pr_i:(V_i,v_i)\rightarrow (U,u_g)$, for all $n'$'s. For each $i$, $pr_i:(V_i,v_i)\rightarrow (U,u_g)$ is of some degree $n_i|n$ and corresponds to some surjective group homomorphism $\pi_1(U, u_g)\xrightarrow[]{\theta_i} \mathbb{Z}/n_i$; fix a $c_i\in\pi_1(U, u_g)$ that maps to $\bar1\in \mathbb{Z}/n_i$ under $\theta_i$. Pick a $pr: (V,v_g)\rightarrow (U,u_g)$ that is a $\mathbb{Z}/{n}$-cover. Suppose it corresponds to $\pi_1(U, u_g)\xrightarrow[]{\theta} \mathbb{Z}/n$ with $c$ the chosen element in $\pi_1(U, u_g)$ above. There is a short exact sequence of groups
\[1\rightarrow\pi_1(V,v_g)\rightarrow\pi_1(U, u_g)\xrightarrow[]{\theta} \mathbb{Z}/n\rightarrow1.\]

Let $Hom(\pi_1(V, v_g),P)$ be the set of group homomorphisms from $\pi_1(V, v_g)$ to $P$. A group homomorphism $\phi\in Hom(\pi_1(V, v_g),P)$ is \textit{$\rho$-liftable}, if there exists a group homomorphism $\widehat{\phi}$ such that the diagram
\[
  \xymatrix @R=.5cm @C=1.3cm {
  \pi_1(V, v_g)\ar@{^{(}->}[d]_{}\ar[r]^{\phi} & P\ar@{^{(}->}[d]\\
  \pi_1(U, u_g)\ar[r]^{\widehat{\phi}} &  G  \ar@{->>}[r]^{Q_P}& \mathbb{Z}/n
 } \eqno{(\ref{glb mdli}.1)}
 \]
\\commutes and the bottom horizontal arrow $\pi_1(U, u_g)\rightarrow\mathbb{Z}/n$ is $\theta$, where $Q_P$ is the projection map. We also say that $\widehat{\phi}$ \textit{lifts} $\phi$. There are two different meanings of ``lift" in the paper (see Remark~\ref{lift}). In this situation, the pointed $P$-cover of $(V,v_g)$ corresponding to $\phi$ is called a \textit{pointed $\rho$-liftable cover of $(V,v_g)$}. If $\widehat{\phi}(c)=(p,\bar 1)$ for some $p\in P$ then $(\phi,p)$ is called a \textit{$\rho$-liftable pair}.

Let $(S,s_0)\in \mathcal{S}_1$. When $(V,U)$ is replaced by $(S\times V,S\times U)$, similarly a \textit{$\rho$-liftable} $\widetilde{\phi}\in Hom(\pi_1(S\times V, (s_0,v_g)),P)$ is defined; the pointed family of $P$-covers of $V$ parameterized by $S$ corresponding to $\widetilde{\phi}$ is called a \textit{pointed $\rho$-liftable family}.

 Denote by $c_{\bullet}$ the image of $c$ under the group homomorphism $\pi_1(U,u_g)\rightarrow\pi_1(S\times U, (s_0,u_g))$ induced by $U\hookrightarrow S\times U$. Similarly a
 \textit{$\rho$-liftable pair} $(\widetilde{\phi},p)$ is defined. A \textit{pointed $\rho$-liftable family pair} means a pair whose first entry is the pointed $\rho$-liftable family corresponding to $\widetilde{\phi}$ and the second entry $p$, for some $(\widetilde{\phi},p)$ a $\rho$-liftable pair. The pair is also denoted by $([\widetilde{\phi}],p)$.\\

Below are two Lemmas for Theorem~\ref{thm3}.

Irreducible linear representations of $\mathbb{Z}/n$ over the field $\mathbb{F}_\mathsf{p}$ correspond to the direct summands in $\mathbb{F}_\mathsf{p}[x]/(x^n-1)=\oplus_i \mathbb{F}_\mathsf{p}[x]/(f_i(x))$, where $f_i(x)$'s are irreducible factors of $x^n-1$ over $\mathbb{F}_\mathsf{p}$. The action of $\bar 1\in \mathbb{Z}/n$ on $\mathbb{F}_\mathsf{p}[x]/(f_i(x))$ is multiplication by $[x]$, where $[x]$ means the equivalence class of $x$. Thus for a pair $(H,\rho)$ in the case of Lemma~\ref{reg}, there is a group isomorphism $H\xrightarrow{\tau}\mathbb{F}_q$, where $q={\mathsf{p}}^m$, such that the induced action of $\rho(-\bar1)$ on $\mathbb{F}_q$ is the multiplication by some $e_\rho\in \mathbb{F}_q$ with $e_\rho^n=1$.

\begin{lem}\label{reg}
Let $P=H=(\mathbb{Z}/\mathsf{p})^m$, an elementary abelian group, and suppose the action $\rho$ on $H$ is irreducible (i.e.\ $\rho$ can not be an action on any subgroup of $H$). A group homomorphism $\phi\in Hom(\pi_1(V, v_g),H)$ is $\rho$-liftable iff for every $b\in\pi_1(V,v_g)$
\[\phi(c^{-1}bc)=\rho(-\bar1)(\phi(b)).\eqno{(*)}\] Moreover, if $\rho=1$, there is only one $\widehat{\phi}$ that can lift $\phi$, and in this case $\widehat{\phi}(c)=(n_{-1}\phi(c^n),\bar1)$, where $n_{-1}$ is a natural number such that $n_{-1}n\equiv1\;(mod\;\textsf{p})$. If $\rho\neq1$, there is a set \{$\widehat{\phi}_h|h\in H$\} consisting of $|H|$ elements that can all lift $\phi$ and in this case $\widehat{\phi}_h(c)=(h,\bar1)$.
\end{lem}
\begin{proof}
\textit{Only if :} If there is a $\widehat{\phi}$ fitting in the diagram of (\ref{glb mdli}.1), then $\phi(c^{-1}bc)=\widehat{\phi}(c)^{-1}\phi(b)\widehat{\phi}(c)$. Since $\widehat{\phi}(c)=(h,\bar1)$ for some $h\in H$, $\phi(c^{-1}bc)=(h,\bar1)^{-1}\phi(b)(h,\bar1)=\rho(-\bar{1})(\phi(b))$.

\textit{If :} Suppose $(*)$ holds. For every element $h\in H$ define a map $\widehat{\phi}_h:\pi_1(U,u_g)\rightarrow H\rtimes_{\rho} \mathbb{Z}/n$ by $\widehat{\phi}_h(bc^i)=\phi(b)(h,\bar1)^i$. The map is well defined since every element in $\pi_1(U,u_g)$ can be written uniquely in the form $bc^i$ with $b\in \pi_1(V,v_g)$ and $0\leq i\leq n-1$. Such $\widehat{\phi}_h$'s are not necessarily homomorphisms; they make the diagram commute. The map $\widehat{\phi}_h$ is a homomorphism iff $\phi(c^n)=(h,\bar 1)^n$. If $\rho=1$, $(h,\bar 1)^n=(nh,\bar 0)$. Then there is a unique $h_0=n_{-1}\phi(c^n)\in H$ such that $\widehat{\phi}_{h_0}$ is a homomorphism. If $\rho\neq1$, the condition automatically holds since both sides equal 0. One can compute $(h,\bar 1)^n=0$ using $\rho(-\bar1)(h)=e_\rho h$. Hence for every $h\in H$, $\widehat{\phi}_h$ is a homomorphism.
\end{proof}

Here is the second lemma needed in the proof of Theorem~\ref{thm3}.

Let $\underline{\sigma}$ be the automorphism in $Gal(V/U)$ corresponding to $\bar1\in \mathbb{Z}/n$. Since $U$ and $V$ are affine, $U=Spec(A)$ and $V=Spec(B)$ for some rings $A$ and $B$. Then $\underline{\sigma}$ corresponds to a ring automorphism $\sigma\in Gal(B/A)$.
\begin{lem}\label{ff}
A group homomorphism $\phi\in Hom(\pi_1(V, v_g),H)$ satisfies condition $(*)$ of Lemma~\ref{reg} iff $\phi$ makes the diagram commutative:
\[
\xymatrix @R=.8cm @C=1.3cm {
 \pi_1(V, v_g)\ar[d]_{\underline{\sigma}_*} \ar[r]^{\phi} & H \ar[d]^{\rho(-\bar{1})} \\
 \pi_1(V, v_{g1}) \ar[r]_{\phi_1} & H
,}\eqno{(\ref{ff}.1)}
\]
where $v_{g1}$ is the image of $v_g$ under $\underline{\sigma}$, $\underline{\sigma}_*$ induced by $\underline{\sigma}$ and $\phi_1$ induced from $\phi$ using any chemin $v_g\rightarrow v_{g1}$.
\end{lem}
\begin{proof}
Since $H$ is abelian, any chemin $v_g\rightarrow v_{g1}$ gives the same isomorphism $\pi_1(V, v_{g1})\simeq \pi_1(V, v_g)$, thus induces the same $\phi_1$ from $\phi$.

Denote by $Fib_{v0}$ (resp.\ $Fib_{v1}$) the fiber functor from (Finite etale covers of $V$) to (Sets) at $v_g$ (resp.\ $v_{g1}$). Similarly denote by $Fib_{u0}$ the fiber functor from (Finite etale covers of $U$) to (Sets) at $u_g$. Denote by $PL_{VU}$ the pullback functor from (Finite etale covers of $U$) to (Finite etale covers of $V$) using $V\rightarrow U$. There are canonical isomorphisms $i_0$ from $Fib_{v0}\circ PL_{VU}$ to $Fib_{u0}$ and $i_1$ from $Fib_{v1}\circ PL_{VU}$ to $Fib_{u0}$.

The element $c\in \pi_1(U, u_g)$ maps to $\bar 1$ under $\theta$, and $\bar1\in \mathbb{Z}/n$ corresponds to $\underline{\sigma}\in Gal(V/U)$, which sends $v_g$ to $v_{g1}$. And since every finite etale cover of $V$ composed with $V\rightarrow U$ is a finite etale cover of $U$, $c\in \pi_1(U, u_g)$ induces $c_{01}$ a chemin $v_{g1}\rightarrow v_g$. So the first and the last squares of diagram (\ref{ff}.2) commute:
\[
\xymatrix @R=.8cm @C=.8cm {
  Fib_{v0}\circ PL_{VU}\ar[d]_{i_0}\ar[r]^{c_{01}}
  & Fib_{v1}\circ PL_{VU}\ar[d]_{i_1}\ar[r]^{\underline{\sigma}_*(b)}
  & Fib_{v1}\circ PL_{VU}\ar[d]^{i_1}\ar[r]^{c_{01}^{-1}}
  & Fib_{v0}\circ PL_{VU}\ar[d]^{i_0} \\
  Fib_{u0}\ar[r]^{c}
  & Fib_{u0}\ar[r]^{\pi_*(b)}
  & Fib_{u0}\ar[r]^{c^{-1}}
  & Fib_{u0}
.}\eqno{(\ref{ff}.2)}
\]
Since $\pi_*(b)=\pi_*\underline{\sigma}_*(b)$ for every $b\in\pi_1(V, v_{g})$,
\[
\xymatrix{
  \pi_1(V, v_{g}) \ar[rr]^{\underline{\sigma}_*} \ar[dr]_{\pi_*}
                &  &    \pi_1(V, v_{g1}) \ar[dl]^{\pi_*}    \\
                & \pi_1(U, u_g)                 .}
\]
The middle square of diagram (\ref{ff}.2) commutes.

Hence the whole diagram (\ref{ff}.2) is commutative, which shows $\phi(c^{-1}bc)=
\phi(c_{01}^{-1}
\underline{\sigma}_*(b)c_{01})$. Since $\phi(c_{01}^{-1}
\underline{\sigma}_*(b)c_{01})
=\phi_1(\underline{\sigma}_*(b))$, the lemma follows.
\end{proof}

The two lemmas above are used to prove Theorem~\ref{thm3}, the first step in the three step construction of the fine moduli space in Theorem~\ref{thm1}.

\begin{dfn}\label{F1}
Define $F_{V,H}^{\rho}$: $\mathcal{S}_1\rightarrow$ (Sets) as the contravariant functor given by $F_{V,H}^{\rho}(S,s_0)$ = \{[$\widetilde{\phi}$]$|\,\widetilde{\phi}:\pi_1(S\times V,(s_0,v_g))\rightarrow H$ is $\rho$-liftable\}, the set of $\rho$-liftable families of $H$-covers of $V$ parameterized by $S$ pointed over $(s_0,v_g)$.

Let $S=Spec(k)$ with $s_0$ determined by $v_g$ using diagram $(\ref{secn&t}.1)$. Then $F_{V,H}^{\rho}(S,s_0)$ is the set of all isomorphism classes of $\rho$-liftable pointed $H$-covers of $(V,v_g)$.
\end{dfn}
\begin{thm}\label{thm3}
Let H be an elementary abelian group $(\mathbb{Z}/\mathsf{p})^m$, $\rho:\mathbb{Z}/n \rightarrow Aut(H)$ an irreducible action of $\mathbb{Z}/n$ on $H$, and $V\rightarrow U$ as above in this section. There is a fine moduli space $M_{V,H}^{\rho}$ representing $F_{V,H}^{\rho}$, the functor for isomorphism classes of pointed $\rho$-liftable $H$-covers of $(V,v_g)$, which is an ind affine space.
\end{thm}
\begin{proof}
In the proof, we will pass between $\mathbb{F}_q$ and $H$ freely using the isomorphism $\tau$ between them given above Lemma~\ref{reg}.

Let $F=F_{V,H}^{\rho}$.

The Artin-Schreier short exact sequence $0\rightarrow\mathbb{F}_q\rightarrow\mathbb{G}_a\xrightarrow[]{\wp}\mathbb{G}_a\rightarrow 0$, where $\wp(f)=f^q-f$, yields $H^0(V,\mathcal{O})\xrightarrow[]{\wp} H^0(V,\mathcal{O})\rightarrow H^1(V,\mathbb{F}_q) \rightarrow 0$, where $0=H^1(V,\mathcal{O})$. This is a short exact sequence of $\mathbb{F}_q$-vector spaces.

Let $\mathbb{X}$ be the subset of $Hom(\pi_1(V, v_g),H)$ that consists of all the isomorphism classes of pointed $\rho$-liftable $H$-covers of $(V,v_g)$. Let $\phi$ be any element in $Hom(\pi_1(V, v_g),H)$. By Lemma~\ref{ff}, $\phi\in \mathbb{X}$ iff $\phi_1\circ\underline{\sigma}_*=e_\rho\phi$. Let $\underline{\sigma}^*: H^1(V,\mathbb{F}_q)\rightarrow H^1(V,\mathbb{F}_q)$ be induced by $\underline{\sigma}$; it is a homomorphism of $\mathbb{F}_q$-vector spaces. Identify $H^1(V,\mathbb{F}_q)$ with $Hom(\pi_1(V, v_g),\mathbb{F}_q)$. By definition of $\underline{\sigma}^*,\underline{\sigma}^*(\phi)=\phi_1\circ\underline{\sigma}_*$.
So $\phi\in \mathbb{X}$ iff
\[\underline{\sigma}^*(\phi)=e_\rho\phi,\eqno{(*1)} \]
which shows that $\mathbb{X}$ is an $\mathbb{F}_q$-subspace of $H^1(V,\mathbb{F}_q)$.

There is a commutative diagram consisting of two short exact sequences of $\mathbb{F}_q$-vector spaces with every symbol already defined above:
\[
\xymatrix @R=.8cm @C=.8cm {
  B=H^0(V,\mathcal{O})\ar[d]_{\sigma}\ar[r]^{\wp}
  & H^0(V,\mathcal{O})\ar[d]_{\sigma}\ar[r]^{\pi}
  & H^1(V,\mathbb{F}_q)\ar[d]^{\underline{\sigma}^*}\ar[r]
  & 0 \\
  H^0(V,\mathcal{O})\ar[r]^{\wp}
  & H^0(V,\mathcal{O})\ar[r]^{\pi}
  & H^1(V,\mathbb{F}_q)\ar[r]
  & 0
,}\eqno{(\ref{thm3}.1)}
\]
which comes from a commutative diagram consisting of two Artin-Schreier short exact sequences of sheaves:
\[
\xymatrix @R=.8cm @C=.8cm {
  0\ar[r]
  &\mathbb{F}_q\ar[d]_{}\ar[r]^{}
  &\mathbb{G}_a\ar[d]_{}\ar[r]^{\wp}
  &\mathbb{G}_a\ar[d]^{}\ar[r]
  & 0 \\
  0\ar[r]
  &\underline{\sigma}_*{F}_q\ar[r]^{}
  &\underline{\sigma}_*\mathbb{G}_a\ar[r]^{\wp}
  &\underline{\sigma}_*\mathbb{G}_a\ar[r]
  & 0
,}
\]
where $\mathbb{F}_q\rightarrow\underline{\sigma}_*
\mathbb{F}_q$ is induced by $V\xrightarrow{\underline{\sigma}} V$ and similarly for $\mathbb{G}_a$.

Let $b\in B$. By the right square of diagram (\ref{thm3}.1), $(\ast1)$ implies
\[\phi:=\pi b\in\mathbb{X}\,\Leftrightarrow\,\sigma b\in e_\rho b+\wp{B}.\eqno{(*2)}\]

Define $D=\sigma-e_\rho$: $B\rightarrow B$, an $A$-module endomorphism of $B$, where $e_\rho$ acts on $B$ by multiplication. Similarly to the proof of Theorem 1.2 in [H80], there is an exact sequence $KerD\xrightarrow[]{\wp}KerD\xrightarrow[]{\pi}\mathbb{X}\rightarrow0$, of $\mathbb{F}_q$-vector spaces. (Denote the restriction of $\wp$ (resp. $\pi$) to $KerD$ also by $\wp$ (resp. $\pi$).)

(\ref{thm3}.1) Now construct $M^{\rho}_{V,H}$ using the $Ker$ short exact sequence above. Let $(KerD)_n=KerD\cap H^0(V,q^n Div_{V})$, where $Div_{V}=\Sigma P_i$ the sum of all the closed points in $\overline{V}- V$ and $\overline{V}$ is the smooth completion of $V$. There is a $k$-vector space filtration $(KerD)_0\leq(KerD)_1\leq ...\leq (KerD)_n\leq...$ . Let $\mathbb{X}_n=\pi((KerD)_n)$. There is a short exact sequence $(KerD)_{n-1}\xrightarrow[]{\wp}(KerD)_n
\xrightarrow[]{\pi}\mathbb{X}_n\rightarrow0$ obtained from the similar one above. Inductively choose bases $ K_n$ of each $(KerD)_n$ as a finite dimensional $k$-vector space, such that $K_{n+1}$ includes both $K_{n}$, and $\{f^q|f\in K_n-K_{n-1}$ and $f$ is not in $k$\}. This is the way to choose bases inductively in a similar situation in the proof of Theorem 1.2 in [H80]. The restriction of $\pi$ to the $k$-linear span $\langle K_{i}-K_{i-1}\rangle_k$ of $K_{i}-K_{i-1}$ is an isomorphism of $\mathbb{F}_q$-vector spaces $\langle K_{i}-K_{i-1}\rangle_k\xrightarrow[]{\pi}\mathbb{X}_i$, which gives a $k$-vector space structure to $\mathbb{X}_i$.

Let $(S,s_0)\in \mathcal{S}_1$ with $S=Spec(R)$. Similarly there is a commutative diagram consisting of two short exact sequences of $\mathbb{F}_q$-vector spaces:
\[
\xymatrix @R=.8cm @C=.8cm {
  H^0(S\times V,\mathcal{O})\ar[d]_{\hat\sigma}\ar[r]^{\wp}
  & H^0(S\times V,\mathcal{O})\ar[d]_{\hat\sigma}\ar[r]^{\Pi}
  & H^1(S\times V,\mathbb{F}_q)\ar[d]^{\hat{\underline{\sigma}}^*}\ar[r] & 0 \\
  H^0(S\times V,\mathcal{O})\ar[r]^{\wp}
  & H^0(S\times V,\mathcal{O})\ar[r]^{\Pi}
  & H^1(S\times V,\mathbb{F}_q)\ar[r]
  & 0
,}
\]
where $\hat\sigma$ is an $R\otimes_k A$-module endomorphism: $R\otimes_k B\rightarrow R\otimes_k B,\, r\otimes b \mapsto r\otimes \sigma(b)$ and $\wp:r\otimes b\mapsto(r\otimes b)^q-r\otimes b$. Let $\widehat{D}=\hat\sigma-e_\rho$. As above, there is a short exact sequence of $\mathbb{F}_q$-vector spaces $Ker\widehat{D}\xrightarrow[]{\wp}Ker\widehat{D}
\xrightarrow[]{\Pi}\widehat{\mathbb{X}}\rightarrow0$, where $\widehat{\mathbb{X}}$ denotes $\{\widetilde\phi\in H^1(S\times V,\mathbb{F}_{q})|\hat{\underline{\sigma}}^*(\widetilde\phi)=
e_\rho\widetilde\phi\}$ and $\hat{\underline{\sigma}}\in Gal(S\times V/S\times U)$ corresponds to $\hat\sigma$. One can check that $Ker\widehat{D}=R\otimes_k Ker D$.

If two $S$-parametrized $P$-covers of $V$ pointed over $(s_0,v_g)$ are equivalent, they are considered the same element in $F(S,s_0)$, by definition of $F$. Hence $F(S,s_0)=\frac{\widehat{\mathbb{X}}+H^1(S,\mathbb{F}_q)}{H^1(S,\mathbb{F}_q)}
=\frac{\widehat{\mathbb{X}}}{\widehat{\mathbb{X}}\cap H^1(S,\mathbb{F}_q)}$ (see Remark~\ref{fml} b). The automorphism $\hat{\underline{\sigma}}$ of $S\times V$ does not change the $S$-factor, thus for any $\tilde{\phi}\in H^1(S,\mathbb{F}_q), \,\, \hat{\underline{\sigma}}^*(\tilde{\phi})=\tilde{\phi}$. If $e_\rho\neq 1,\, \widehat{\mathbb{X}}\cap H^1(S,\mathbb{F}_q)=0$ and $F(S,s_0)=\widehat{\mathbb{X}}=\Pi( Ker\widehat{D})=\Pi(R \otimes_k Ker D) $. Let the transition map from $R\otimes_k \langle K_n-K_{n-1}\rangle_k$ to $R\otimes_k \langle K_{n+1}-K_n\rangle_k$ be Frobenius ($r\otimes b\mapsto (r\otimes b)^q$). Then $\varinjlim_n R\otimes_k \langle K_n-K_{n-1}\rangle_k$ is an $\mathbb{F}_q$-vector space. There is an $\mathbb{F}_q$-vector space isomorphism $\varinjlim_n R\otimes_k \langle K_n-K_{n-1}\rangle_k\xrightarrow{\Pi}\Pi(R \otimes_k Ker D) $. Hence $F(S,s_0)=\varinjlim_n R\otimes_k \langle K_n-K_{n-1}\rangle_k$. Write out elements in $ K_n-K_{n-1}$ as $\{\verb"k"_1,...,\verb"k"_{d_K}\}$, whose dual vectors are $\{\verb"k"_1\spcheck,...,
\verb"k"_{d_K}\spcheck\}$. Define $\mathbb{A}(\langle K_n-K_{n-1}\rangle_k)
$ as $Spec(k[\verb"k"_1\spcheck,...,
\verb"k"_{d_K}\spcheck])$. Now $R\otimes_k \langle K_n-K_{n-1}\rangle_k =Hom_k(\langle K_n-K_{n-1}\rangle_k\spcheck,R)=
Hom_k(S,\mathbb{A}(\langle K_n-K_{n-1}\rangle_k))$. Therefore $M^{\rho}_{V,H}:=$ the ind scheme $\{\mathbb{A}(\langle K_n-K_{n-1}\rangle_k)\}$, where the transition morphism between $\mathbb{A}(\langle K_n-K_{n-1}\rangle_k)$ and $\mathbb{A}(\langle K_{n+1}-K_n\rangle_k)$ is given by Frobenius, represents $F$.

If $e_\rho=1$, then $F(S,s_0)=\frac{H^1(S\times U,\mathbb{F}_q)}{H^1(S,\mathbb{F}_q)}$. This is the case, if $H=\mathbb{Z}/\textsf{p}$, of the base step in the proof of Theorem 1.2 in [H80]; the proof there also works for any elementary abelian group $H$. Hence $F$ is represented by $M_{V,H}^{\rho}:=M_{U,H}$, which is denoted by $M_G$ there with $G=H$, an ind affine space with transition morphisms given by Frobenius as well. Since now $G$ is a product $H\times \mathbb{Z}/n$, it can be derived directly that $F$ is represented by $M_{U,H}$,
\end{proof}

\begin{rem}\label{clpt}
a. By Theorem~\ref{thm3}, for any pointed affine connected $k$-scheme $(S,s_0)$, there is a bijection between $F(S,s_0)$ and $M^{\rho}_{V,H}(S)$, where the latter set is the set of $k$-morphisms from $S$ to $M^{\rho}_{V,H}$.

b. Let $S=Spec(k)$ with $s_0$ determined by $v_g$ using diagram $(\ref{secn&t}.1)$. Then $F(S,s_0)$, the set of all $\rho$-liftable pointed $H$-covers of $(V,v_g)$, are in bijection with $M_{V,H}^{\rho}(S)$, the set of $k$-points of $M_{V,H}^{\rho}$, same as the set of closed points of $M_{V,H}^{\rho}$.

c. Let $M^{\rho}_{V,H,n}=Spec(k[\verb"k"_1
\spcheck,...,
\verb"k"_{d_K}\spcheck])$ if $\rho\neq1$. It is the \textit{$n$-th piece} of $M^{\rho}_{V,H}$. Similarly the $n$-th piece of $M^{\rho}_{V,H}$ when $\rho=1$ can be defined.
\end{rem}

\begin{rem}\label{uni fml}
(Remark/Definition)

With the same notations of Theorem~\ref{thm3}.

Let $M^{\rho}_{V,H,n}$ be the $n$-th piece of $M^{\rho}_{V,H}$ (see Remark~\ref{clpt} c). \textit{A compatible system of $H$-covers of $V$ over $M^{\rho}_{V,H}$} means a collection of covers \{$H$-covers $\widetilde{Z}_n$ of $M^{\rho}_{V,H,n}\times V\,|\,n\geq1$\} such that $\widetilde{Z}_n$ pulled back to $M^{\rho}_{V,H,n-1}\times V$ is isomorphic to $\widetilde{Z}_{n-1}$. Since $H$ is abelian, given any point $m$ on $M^{\rho}_{V,H,n}$, where we point $\widetilde{Z}_n$ over $(m,v_g)$ does not matter by Remark~\ref{abel}.

\textit{A universal family representative over the moduli space} $M^{\rho}_{V,H}$ means, a compatible system of $H$-covers \{$H$-covers $\widetilde{Z}_n$ of $M^{\rho}_{V,H,n}\times V\,|\,n\geq1$\} of $V$ over $M^{\rho}_{V,H}$, which can be used to give the isomorphism of functors $Hom(\bullet,M^{\rho}_{V,H})
\xrightarrow{\simeq}F_{V,H}^{\rho}
:\mathcal{S}_1\rightarrow(Sets)$ given in Theorem~\ref{thm3} as follows: Sending a $k$-morphism $f$ from a $k$-scheme $S$ with $(S,s_0)\in \mathcal{S}_1$ to $M^{\rho}_{V,H}$ (see Definition~\ref{indschmorphdfn} a), to the equivalence class of the pullback of $\widetilde{Z}_n$, using the morphism $f$, to $S\times V$ pointed arbitrarily over $(s_0,v_g)$, is the isomorphism of functors $Hom(\bullet,M^{\rho}_{V,H})
\xrightarrow{\simeq}F_{V,H}^{\rho}$ given in Theorem~\ref{thm3}. Since the $\widetilde{Z}_n$'s are compatible any $n$ can be used.

It is derived from definitions that any two universal family representatives are equivalent. The equivalence class of a universal family representative is \textit{the universal family over the moduli space $M^{\rho}_{V,H}$}. There must be a universal family representative: Let $S$ be $M^{\rho}_{V,H,n}$. Identity morphism of $M^{\rho}_{V,H,n}$ determines a morphism $S\rightarrow M^{\rho}_{V,H}$. The morphism gives an equivalence class in $F_{V,H}^{\rho}(S,m)$ using $Hom(\bullet,M^{\rho}_{V,H})
\xrightarrow
{\simeq}F_{V,H}^{\rho}$ given in Theorem~\ref{thm3}, for any point $m$ on $S$. Then $\rho$-liftable representatives in the equivalence class are candidates for the $n$-th element of a universal family representative. Use the same kind of argument as in Lemma 4.25 of [TY17], a compatible system of $H$-covers can be chosen.

If $\rho\neq1$, a universal family representative over the moduli space $M^{\rho}_{V,H}$ can be given by \{the $H$-cover of $M^{\rho}_{V,H,n}\times V$ given by $z^q-z=\sum_{k_i\in K_{n}-K_{n-1}} k_i\spcheck\otimes k_i$$\,|\,n\geq1$\}, by the construction of $M^{\rho}_{V,H}$. The $H$-covers are compatible for different $n$'s.

If $\rho=1$, similarly a universal family representative over $M^{\rho}_{V,H}$ can be given explicitly: Replace $z^q-z=\sum_{k_i\in K_{n}-K_{n-1}} k_i\spcheck\otimes k_i$ above by $z^q-z=\sum_{l_i\in L_{n}-L_{n-1}} l_i\spcheck\otimes l_i$. $L_{n}$, an analogue of $ K_{n}$, is the basis chosen inductively for $A_n/{k^+}$ in the proof of Theorem 1.2 in [H80]; here $A_n=H^0(U,q^n Div_{U})$ with $\textsf{p}$ there replaced by $q$ and $B_n$ there is denoted by $L_n$ here.

The universal family representative over $M^{\rho}_{V,H}$ given above is \textit{the canonical universal family representative over $M^{\rho}_{V,H}$}. The $n$-th $H$-cover in every other universal family representative over $M^{\rho}_{V,H}$, differs from the $n$-th $H$-cover in the canonical one by an element in $H^1(M^{\rho}_{V,H,n},H)$, as shown in the last two paragraphs in the proof of Theorem~\ref{thm3}.
\end{rem}

The corollary below is a version of Theorem~\ref{thm3} for pairs, which will be used in the proof of Theorem~\ref{thm2}.

\begin{dfn}\label{F1pair}
With the same setting as in Theorem~\ref{thm3}. Define $F_{V,H}^{\rho\bullet}$: $\mathcal{S}_1\rightarrow$ (Sets) as the contravariant functor given by $F_{V,H}^{\rho\bullet}(S,s_0)$ = \{([$\widetilde{\phi}$],$h$)
$|\,\widetilde{\phi}:\pi_1(S\times V,(s_0,v_g))\rightarrow H$ and $(\widetilde{\phi},h)$ is a $\rho$-liftable pair\}, the set of $\rho$-liftable family pairs of $H$-covers of $V$ parameterized by $S$, pointed over $(s_0,v_g)$.

Let $S=Spec(k)$ with $s_0$ determined by $v_g$ using diagram $(\ref{secn&t}.1)$. Then $F_{V,H}^{\rho\bullet}(S,s_0)$ is the set of $\rho$-liftable pairs of $V$.
\end{dfn}
\begin{cor}\label{disj uni over H}
Under the same setting of Theorem~\ref{thm3}, there is a fine moduli space $M_{V,H}^{\rho\bullet}$ representing $F_{V,H}^{\rho\bullet}$, the functor for $\rho$-liftable pairs of $V$. It is a disjoint union of finitely many copies of $M_{V,H}^{\rho}$ in Theorem~\ref{thm3}.
\end{cor}
\begin{proof}
Let $M_{V,H}^{\rho\bullet}$ be $M_{V,H}^{\rho}$ if $\rho=1$ and $\amalg_{h\in H} M_{V,H,h}^{\rho}$ if $\rho\neq1$, where $M_{V,H,h}^{\rho}$ means a copy of $M_{V,H}^{\rho}$ indexed by $h$.

Let $\phi\in Hom(\pi_1(V,v_g),H)$ and $(\phi,h_0)$ be a $\rho$-liftable pair. The map $\phi_{h_0}$, as defined in the proof of Lemma~\ref{reg}, is in fact a homomorphism. As stated in Lemma~\ref{reg}, if $\rho=1$ then $h_0$ is the only element in $H$ such that $(\phi,h_0)$ is a $\rho$-liftable pair. If $\rho\neq 1$ then for every $h\in H$ the pair $(\phi,h)$ is $\rho$-liftable. Using this fact and Theorem~\ref{thm3} $F_{V,H}^{\rho\bullet}$ is represented by $M_{V,H}^{\rho\bullet}$.
\end{proof}
By Corollary~\ref{disj uni over H}, $M_{V,H}^{\rho}$ is a connected component of the ind scheme $M_{V,H}^{\rho\bullet}$. See Remark~\ref{compo}.

Here is the 2nd step of the 3 step construction of the fine moduli space in Theorem~\ref{thm1}. Let $P$ be an arbitrary finite $\mathsf{p}$-group now.
\begin{dfn}\label{F2}
Let $F_{V,P}^{\rho\bullet}$: $\mathcal{S}_1\rightarrow$ (Sets) be the contravariant functor given by $F_{V,P}^{\rho\bullet}(S,s_0)$ = \{([$\widetilde{\phi}$],$p$)
$|\,\widetilde{\phi}:\pi_1(S\times V,(s_0,v_g))\rightarrow P$ and $(\widetilde{\phi},p)$ is a $\rho$-liftable pair\}, the set of $\rho$-liftable family pairs of $P$-covers of $V$ parameterized by $S$, pointed over $(s_0,v_g)$.

Let $S=Spec(k)$ with $s_0$ determined by $v_g$ using diagram $(\ref{secn&t}.1)$. Then $F_{V,P}^{\rho\bullet}(S,s_0)$ is the set of $\rho$-liftable pairs of $V$.
\end{dfn}

\begin{dfn}\label{unifmlpair}
A similar definition for pairs to a universal family representative over $M^{\rho}_{V,H}$ in Remark~\ref{uni fml} will be given.

Assume there is an ind scheme $M^{\rho\bullet}_{V,P}$, consisting of finitely many connected components (see Remark~\ref{compo}), representing $F_{V,P}^{\rho\bullet}$
with an isomorphism between functors $Hom(\bullet,M^{\rho\bullet}_{V,P})
\xrightarrow{\simeq} F_{V,P}^{\rho\bullet}$.

Below $M^{\rho\bullet}_{V,P}$ is viewed as a scheme instead of an ind scheme (see Remark~\ref{n}). Connected components of $M^{\rho\bullet}_{V,P}$ are denoted by $\{M^{\rho}_{V,P,j}\}$.

\textit{A system of universal family pair representatives over} $M^{\rho\bullet}_{V,P}$, means a collection of a $\rho$-liftable pair $(\tilde\phi_{0j,m_j}:\pi_1(M^{\rho}_{V,P,j}\times V,(m_j,v_g))\rightarrow P,p_{\tilde\phi_{0j,m_j}})$ for every base point $m_j$ over each $M^{\rho}_{V,P,j}$, which can be used to give the isomorphism of functors $Hom(\bullet,M^{\rho\bullet}_{V,P})
\xrightarrow{\simeq}
F^{\rho\bullet}_{V,P}$ as follows: Sending a $k$-morphism $S\xrightarrow{\mathfrak{c}} M^{\rho}_{V,P,j}$ with $(S,s_0)\in \mathcal{S}_1$ and $s_0$ mapped to $m_j$ under $\mathfrak{c}$, to the pair ([$\tilde\phi_{0j,m_j}\circ \tilde{\mathfrak{c}}_*$],$p_{\tilde\phi_{0j,m_j}}$) with $S\times V\xrightarrow{\tilde{\mathfrak{c}}}M^{\rho}_{V,P,j}\times V$ induced by $\mathfrak{c}$ and $\tilde{\mathfrak{c}}_*$ the homomorphism between fundamental groups induced by $\tilde{\mathfrak{c}}$, is the given isomorphism of functors $Hom(\bullet,M^{\rho\bullet}_{V,P})
\xrightarrow{\simeq}
F_{V,P}^{\rho\bullet}$.

It can be derived from definition that any two $\tilde\phi_{0j,m_j}$ and $\tilde\phi'_{0j,m_j}$ in two different systems of universal family pair representatives over $M^{\rho\bullet}_{V,P}$ are equivalent and $p_{\tilde\phi_{0j,m_j}}=p_{\tilde\phi'_{0j,m_j}}$.

Similarly to Remark~\ref{uni fml}, there must be a system of universal family pair representatives over $M^{\rho\bullet}_{V,P}$.

\end{dfn}
\begin{thm}\label{thm2}
With the notations above, there exists a fine moduli space $M_{V,P}^{\rho\bullet}$ representing $F_{V,P}^{\rho\bullet}$, the functor for $\rho$-liftable pairs of $V$. It is a disjoint union of finitely many ind affine spaces.
\end{thm}
\begin{proof}
Induct on $|P|$.

Let $F=F_{V,P}^{\rho\bullet}$.

(\ref{thm2}.1) Take a minimal normal subgroup $H$ of $G$ inside $C(P)$,
the nontrivial center of $P$. It is a product of copies of some
simple group $S$. Hence $H\approx(\mathbb{Z}/p)^{m}$ for
some $m\geq1$. Let $\rho_0: \mathbb{Z}/n \rightarrow
Aut(H)$ be the $\mathbb{Z}/n$-action induced by
$\rho$; $\rho_0$ is irreducible by the minimality of $H$.
If $H=P$, then this is in the case of
Corollary~\ref{disj uni over H}. Hence assume $H<P$ below.
Let $\bar \rho: \mathbb{Z}/n \rightarrow Aut(\bar P)$ be the $\mathbb{Z}/n$-action induced by $\rho$, where $\bar P=P/H$. By the inductive hypothesis and Corollary~\ref{disj uni over H} respectively $\bar M:=M^{\bar \rho\bullet}_{V,\bar P}$ and $M^0:=M^{\rho_0\bullet}_{V,H}$ exist.

Denote $F^{\bar \rho\bullet}_{V,\bar P}$ by $\bar F$.

It will be shown that $\bar M\times M^0$ is the moduli space
desired.

(\ref{thm2}.2) First need to lift a system of universal family pair
representatives over $\bar M$.  For every point $\bar m$ of
$\bar M$, there is a $\bar \rho$-liftable pair ($
\widetilde{\mu}_{0}: \pi_1(\bar M\times V,(\bar m,v_g))
\rightarrow \bar P,\bar{p}_0)$ in the system, which is the
counterpart of $(\tilde\phi_{0j,m_j}:\pi_1(M^{\rho}_{V,P,j}
\times V,(m_j,v_g))\rightarrow \bar P,p_{\tilde\phi_{0j,m_j}})$
in Definition~\ref{unifmlpair}. Denote by $c_{\bullet}$ the image of $c$ under the group homomorphism $\pi_1(U,u_g)\rightarrow\pi_1(\bar{M}\times U, (\bar{m},u_g))$ induced by $U\hookrightarrow \bar{M}\times U$. The pair gives a $\widehat{
\widetilde{\mu}}_{0}: \pi_1(\bar M\times U,(\bar m,u_g))
\rightarrow \bar P\rtimes_{\bar \rho} \mathbb{Z}/n$ with
$\widehat{\widetilde{\mu}}_{0}(c_{\bullet})=(\bar{p}_0,\bar1)$,
similar to the diagram (\ref{glb mdli}.1). As $\pi_1(\bar M\times
U,(\bar m,u_g))$ has $cd_p\leq1$ ([H80], p1101), $\widehat{
\widetilde{\mu}}_{0}$ lifts (a different meaning of ``lift", see Remark~\ref{lift}) to a $\widehat{\widetilde{\psi}}_{0}:\pi_1(\bar M\times U,(\bar m,u_g))\rightarrow P\rtimes_{ \rho} \mathbb{Z}/n$ ([Serre], I Prop. 16) with $\widehat{\widetilde{\psi}}_{0}(c_{\bullet})=(p_0,\bar1)$, for some $p_0\in P$ mapping to $\bar p_0\in \bar P$, under the quotient map $P\twoheadrightarrow \bar P$. Denote the restriction of $\widehat{\widetilde{\psi}}_{0}$ on $\pi_1(\bar M\times V,(\bar m,v_g))$ by $\widetilde{\psi}_{0}$.

(\ref{thm2}.3) Then use the lift $\widetilde{\psi}_{0}$ obtained above to
separate a $\rho$-liftable pair of $S\times V$ into two parts.
Let $(S,s_0)\in\mathcal{S}_1$ and suppose ($\widetilde \phi:
\pi_1(S\times V,(s_0,v_g))\rightarrow P$, $p_1$) is a $\rho$-liftable pair. Its quotient ($\bar{ \widetilde {\phi}}:\pi_1(S
\times V,(s_0,v_g))\rightarrow \bar P,\bar p_1)$ is a $\bar
\rho$-liftable pair. By the inductive hypothesis, the quotient pair
corresponds to a morphism $\beta:S\rightarrow \bar M$. Denote $\beta\times Id_{V}$ by $\tilde{\beta}$. Denote
 the induced homomorphism $\pi_1(S\times V,(s_0,v_g))
 \rightarrow \pi_1(\bar M\times V,(\beta(s_0),v_g))$ by ${
 \tilde\beta}_*$, and let $\widetilde{\psi}:=\widetilde{\psi}_0
 \circ {\tilde\beta}_*$ (letting $\bar m$ above be $\beta(s_0)$
 here). Then define a ``quotient homomorphism"
$\widetilde{\eta}: \pi_1(S\times V,(s_0,v_g))\rightarrow H$ by
$\widetilde \phi\widetilde{\psi}^{-1}$: Since $\bar M$ is a fine
moduli space for $\bar F$ that involves equivalence classes,
$\bar{ \widetilde {\phi}}$ and $\bar{\widetilde{\psi}}$ only
agree pulled back to some finite etale cover $T$ of $S$, by
definition of $\bar F$. Pick a point $t_0$ on $T$ mapping to
$s_0$. Define $\widetilde{\eta}_T(a)=\widetilde \phi_T(a)
\widetilde{\psi}_T(a^{-1})$ for every $a\in \pi_1(T\times V,
(t_0,v_g))$, where $\widetilde \phi_T$ means $\widetilde \phi$
 pulled back to $T$ and similarly for $\widetilde{\psi}_T$.
 Actually $\widetilde{\eta}_T$ maps to $H$ and the centrality of
  $H$ in $P$ implies that $\widetilde{\eta}_T$ is a homomorphism
  . Let $h_1=p_1 p_0^{-1}$. Then $(\widetilde{\eta}_T,h_1)$ is a
  $\rho_0$-liftable pair and hence corresponds to a morphism
  $\alpha_T:T\rightarrow M^0$. By etale decent ([H80], p1109,
  second paragraph) $\alpha_T$ descends to a morphism $
  \alpha:S\rightarrow M^0$. Hence get $(\alpha,\beta):S
  \rightarrow M$, where $M=M_{V,P}^{\rho\bullet}:=M^0\times
  \bar M$.

It is straightforward to verify that the assignment
($\widetilde \phi$, $p_1$)$\leadsto$ $(\alpha,\beta)$
is well defined on $\rho$-liftable family pairs (i.e.\
is independent of the choice of $\widetilde \phi$ in its
equivalence class), and yields a bijection between $F(S,s_0)$ and
$Hom(S,M)$. As the bijection is compatible with pullback, it
follows that $M$ represents $F$.
\end{proof}

\begin{rem}\label{unifml2}
Keeping track of universal family representatives in the inductive construction,
a $\rho$-liftable universal family representative
over (each connected component of) $M_{V,P}^{\rho\bullet}=M^
0\times \bar M$ can be given by the product of a $\rho_0$-liftable universal family representative  over $M^0$ and a $
\rho$-liftable lift of a $\bar{\rho}$-liftable universal family
representative over  $\bar M$, using the inclusions $M^0
\hookrightarrow M^0\times \bar M$ and $\bar M
\hookrightarrow M^0\times \bar M$.
\end{rem}
\begin{rem}\label{n}
Since the moduli spaces are ind schemes, strictly speaking
the argument above needs to be carried out for each $n$ and
check compatibility for different $n$'s. The argument given
above has the advantage of being more concise, which follows
the way of presentation in Theorem 1.2 of [H80].
\end{rem}
\begin{rem}\label{dcps}
In the inductive proof of Theorem~\ref{thm2},
$(P,\rho)$ is said to \textit{be decomposed} to $(H,\rho_0)$
and $(\bar P,\bar \rho)$. If $(\bar P,\bar \rho)$ is not in the
case considered in Theorem~\ref{thm3}, then it can be further
decomposed similarly. Repeat the inductive step in Theorem~
\ref{thm2} until the last pair got is in the case of Theorem~
\ref{thm3}. The pairs got in the process are denoted by
$(H_t,\rho_t)_t$, which are all in the case of Theorem~
\ref{thm3} and the first of which is $(H,\rho_0)$. Then
$M_{V,P}^{\rho\bullet}=\Pi_t M_{V,H_t}^{\rho_t\bullet}$,
which consists of finitely many connected components of the form
$\Pi_t M_{V,H_t}^{\rho_t}$.
\end{rem}

Here is the main theorem of this section, on moduli for covers with a given cyclic-by-$\mathsf{p}$ Galois group.
\begin{dfn}\label{F3}
Let $F_{U,G}$: $\mathcal{S}_1\rightarrow$ (Sets) be the contravariant functor given by $F_{U,G}(S,s_0)$ = \{$[\tilde{\phi}]\,|\,\tilde{\phi}:\pi_1(S\times U,(s_0,u_g))\rightarrow G$\}, the set of families of $G$-covers of $U$ parametrized by $S$, pointed over $(s_0,u_g)$.

Let $S=Spec(k)$ with $s_0$ determined by $u_g$ using diagram $(\ref{secn&t}.1)$. Then $F_{U,G}(S,s_0)$ is the set of pointed $G$-covers of $(U,u_g)$.
\end{dfn}
\begin{thm}\label{thm1}
There exists a fine moduli space $M_{U,G}$ representing $F_{U,G}$, the functor for pointed $G$-covers of $(U,u_g)$, which is a disjoint union of finitely many ind affine spaces.
\end{thm}
\begin{proof}
It will be shown that $F_{U,G}$ is isomorphic to $\amalg_{V_i}F_{V_i,P}^{\rho_{n_i},\bullet}$ (see Section~\ref{table of symbols} ``Table of symbols" for $V_i$).  The disjoint union of functors means taking disjoint union of sets, since the functors map to the category of sets. Hence it is represented by $\amalg_{V_i}M_{V_i,P}^{\rho_{n_i},\bullet}$, by Theorem~\ref{thm2}.

First the left to right direction map is given in the isomorphism
wanted.

Let $(S,s_0)\in\mathcal{S}_1$ and $(\widetilde W,\widetilde{w}
_g) \rightarrow (S\times U,(s_0,u_g))$ be a pointed $G$-cover
corresponding to some $\widetilde{\widetilde \phi}\in Hom(\pi_1
(S\times U,(s_0,u_g)),G)$. Let $(\widetilde W_m,\overline{
\widetilde{w}_g})$ be the pointed connected component of $
\widetilde W/P$, a $(\mathbb{Z}/n')$-cover of $(S\times U,(
s_0,u_g))$ with $(\mathbb{Z}/n')$ the order $n'$ subgroup in
$\mathbb{Z}/n$ for some $n'|n$. The diagram commutes:
\[
\xymatrix @R=.5cm @C=1.8cm {
  \pi_1(\widetilde W_m,\overline{\widetilde{w}_g})
  \ar@{^{(}->}[d]_{}\ar[r]^{\widetilde{\phi}_m} & P\ar@{^{(}->}[d]\\
  \pi_1(S\times U,(s_0,u_g))\ar[r]^{\widetilde{\widetilde{\phi}}} & G
.}
\]

Let $T$ be a connected component of the inverse image in
$\widetilde W_m$ of $S\times\{u_g'\}$, where $u_g'$ is any
$k$-point on $U$. The fibers of $\widetilde W_m$ over $k$-points of $U$ does not vary since the degree of the cover is
prime to $\textsf{p}$. The $k$-scheme $T$ is a finite etale
cover of $S$ and pick any base point $t_0$ that maps to
$s_0$. The cover $\widetilde W_{m}$ pulled back to $T\times
U$ is isomorphic to a disjoint union of copies of a product $T
\times V_i$ for some $V_i$ a $\mathbb{Z}/n_i$-cover of $U$:
\[
\xymatrix @R=.5cm @C=1.3cm {
   \coprod (T\times V_i,(t_0,v_i))\ar[d]\ar[r]& (\widetilde W_{m},\overline{\widetilde{w}_g})\ar[d]\\
  (T\times U,(t_0,u_g))\ar[r]& (S\times U,(s_0,u_g))
  ,}
\]
as $(\mathbb{Z}/n')$-covers of $T\times U$, using the canonical
embedding $\iota_{n_i}$ of $\mathbb{Z}/n_i$ in $\mathbb{Z}
/n$ given in Section~\ref{secn&t} Definition~\ref{gpextn} b.

Let $\widetilde{\phi}_T$ be the composition $\pi_1(T\times V_i,
(t_0,v_i))\rightarrow\pi_1(\widetilde W_{m},\overline{\widetilde
{w}_g})\rightarrow P$ induced by $T\times V_i\rightarrow
\widetilde W_{m}$. Let $c_{i\bullet}$ be the image of $c_i$
under $\pi_1(U,u_g)\rightarrow \pi_1(S\times U,(s_0,u_g))$.
Let $p_0$ be the first entry of $\widetilde{\widetilde \phi}(c_
{i\bullet})\in G=P\rtimes_{\rho} \mathbb{Z}/n$. Then $(
\widetilde{\phi}_T,p_0)$ is a $\rho_{n_i}$-liftable pair. It
corresponds to a morphism $\mathfrak{c}_T: T\rightarrow M_
{V_i,P}^{\rho_{n_i},\bullet}$. By etale decent again $\mathfrak{c}
_T$ descends to a morphism $\mathfrak{c}_S:S\rightarrow M_
{V_i,P}^{\rho_{n_i},\bullet}$. The morphism $\mathfrak{c}_S$
corresponds to an element in $F_{V_i,P}^{\rho_{n_i},\bullet}(S,s_
0)$. In fact a morphism $\delta: F_{U,G}\rightarrow
\amalg_{V_i}F_{V_i,P}^{\rho_{n_i},\bullet}$ is got.

Conversely, suppose $(\widetilde{\phi},p_0)$ is a $\rho_{n_i}$-
liftable pair with $\pi_1(S\times V_i,(s_0,v_i))\xrightarrow{
\widetilde{\phi}} P$. The diagram commutes:
\[
\xymatrix @R=.5cm @C=1.3cm {
  \pi_1(S\times V_i,(s_0,v_i))\ar@{^{(}->}[d]_{}\ar[r]^{\widetilde{\phi}} & P\ar@{^{(}->}[d]\\
  \pi_1(S\times U,(s_0,u_g))\ar[r]^{\widehat{\widetilde{\phi}}} & P\rtimes_{\rho_{n_i}} \mathbb{Z}/n_i\lhook\joinrel\xrightarrow{\widetilde{\iota_{n_i}}}G
,}
\]
where $\widehat{\widetilde{\phi}}$ sends $c_{i\bullet}$ to $(p_
0,\bar1)$, and $\widetilde{\iota_{n_i}}$ is the group embedding induced
by $\iota_{n_i}$. Hence a pointed family of $G$-covers of $
U$ parametrized by $S$ corresponding to $\widetilde{\iota_{n_
i}}\circ\widehat{\widetilde{\phi}}$ is got. In fact a morphism $
\gamma: \amalg_{V_i}F_{V_i,P}^{\rho_{n_i},\bullet} \rightarrow
F_{U,G}$ is got, which is inverse to $\delta$.
\end{proof}

\section{Moduli for $\mathsf{p}'$-by-$\mathsf{p}$ covers} \label{p' mdli}
In Section~\ref{p' mdli}, it is shown that, given a pointed affine curve $(U,u_g)$, an intersection of finitely many fine moduli spaces for cyclic-by-$\mathsf{p}$ covers of some affine curves gives a moduli space for $\mathsf{p}'$-by-$\mathsf{p}$ covers of the curve (Corollary~\ref{uniirr}).

The next simplest groups after cyclic-by-$\mathsf{p}$ groups are $\mathsf{p}'$-by-$\mathsf{p}$ groups. The first idea on how to get a moduli space for $\mathsf{p}'$-by-$\mathsf{p}$ covers of $(U,u_g)$, out of fine moduli spaces for cyclic-by-$\mathsf{p}$ covers of affine curves constructed in Section~\ref{glb mdli}, is to intersect them.

The fine moduli spaces for cyclic-by-$\mathsf{p}$ covers of some affine curves, intersect in a fixed fine moduli space $M_{V',P,0}$ for some affine curve $V'$, which is given first below in Remark~\ref{M's}.

Lemma~\ref{(M)} and Lemma~\ref{cmap} show how to embed a fine moduli space for cyclic-by-$\mathsf{p}$ covers of an affine curve in $M_{V',P,0}$. The first lemma is the base case for the induction in the proof of the 2nd lemma.

Then an intersection in $M_{V',P,0}$ gives a target moduli space $M^{0\rho'}_{V',P}$ in Definition~\ref{M}. However, it is not a moduli space for covers of $(U,u_g)$ with Galois group the $\mathsf{p}'$-by-$\mathsf{p}$ group given, because pieces do not patch together well when $P$ is not abelian (see Remark~\ref{gal gp}). It is a moduli space for something else; see Proposition~\ref{intersect mdli}. Similarly pieces may not patch together well for a disconnected $P$-cover. Therefore $M^{0\rho'}_{V',P}$ only contains connected covers. The moduli space for covers of $(U,u_g)$ with Galois group the given $\mathsf{p}'$-by-$\mathsf{p}$ group is a corollary of Proposition~\ref{intersect mdli}.

One final thing for the intersection idea to work, is to use a weaker definition of equivalence. A new $ER$-equivalence is introduced below in the definition of $F^{er,Gal/\rho'}_{V',P}$, the functor to present, and that of $M^{er0\rho'}_{V',P}$, a functor related to the moduli space $M^{0\rho'}_{V',P}$. Using $ER$-equivalence $F^{er,Gal/\rho'}_{V',P}$ and $M^{er0\rho'}_{V',P}$ are proven isomorphic in Proposition~\ref{intersect mdli}.

As always, we follow notations and terminology defined in Section~\ref{secn&t}.\\

First the space where intersections take place is given.

Let $(V',v'_g)\rightarrow (U,u_g)$ be a pointed connected $P'$-cover of $(U,u_g)$, which corresponds to a surjective group homomorphism $\theta': \pi_1(U, u_g)\rightarrow P'$.
\begin{rem}\label{M's}
Since $P$ can be decomposed in different ways in the construction of $M_{V',P}$ (see proof of Theorem 1.2 in [H80]; see Remark~\ref{dcps} for a $\rho$-liftable version), there are different forms of $M_{V',P}$. Since they are all fine moduli spaces of $F_{V',P}$, it is derived from the definition that they are isomorphic. Fix a fine moduli space $M_{V',P,0}$ for $F_{V',P}$ below, where intersections take place.
\end{rem}

Now the objects that intersect later are given.

Let $(V'_{i'},\overline{v'_{gi}})$ be the quotient of $(V',v'_g)$ by $\langle p_i'\rangle$, the subgroup generated by $p_i'$, and let $\rho_i':\langle p_i'\rangle \rightarrow Aut(P)$ be the restriction of $\rho'$. There is a short exact sequence of groups:
\[1\rightarrow\pi_1(V',v'_g)
\rightarrow\pi_1(V'_i, \overline{v'_{gi}})\xrightarrow[]
{\theta'_i} \langle p_i'\rangle\rightarrow1.\]
Let $\pi_{i*}'$ be the homomorphism between fundamental groups induced by $\pi_i':V'_i\rightarrow U$. The following diagram commutes:
\[\xymatrix{
 \pi_1(V'_i,\overline{v'_{gi}})  \ar[d]_{\pi_{i*}'} \ar[r]^{\theta'_i} & \langle p_i'\rangle\ar[d]^{\subset} \\
  \pi_1(U,u_g)\ar[r]^{\theta'} & P'   .}
\]
For every $p_i'\in P'$, fix a $c'_i$ in $\pi_1(V'_i,\overline{v'_{gi}})$ that maps to $p_i'$ under  $\theta_i'$. The pointed $\langle p_i'\rangle$-cover $(V',v'_g)\rightarrow (V'_{i'},\overline{v'_{gi}})$ is the counterpart of the pointed $\mathbb{Z}/n$-cover $(V,v_g)\rightarrow (U,u_g)$ in Theorem~\ref{thm2} of Section~\ref{glb mdli}. Apply Theorem~\ref{thm2} on $(V',v'_g)\rightarrow (V'_{i'},\overline{v'_{gi}})$ and a fine moduli space $M^{\rho'_i,\bullet}_{V',P}$ for $\rho'_i$-liftable pairs of $(V',v'_g)$ is got.

For every $p'_i$ denote by $\{M^{\rho'_i}_{V',P,ij}\}$ the set of finitely many connected components of $M^{\rho'_i,\bullet}_{V',P}$. Denote by $(M^{\rho'_i}_{V',P,ij})_i$ a tuple of connected components indexed by $i$, an element in $\Pi_i\{M^{\rho'_i}_{V',P,ij}\}$. For each tuple $(M^{\rho'_i}_{V',P,ij})_i$ do their intersection in $M_{V',P,0}$, the way of which will be defined below. Then take the disjoint union of intersections belonging to different tuples. The disjoint union is almost $M^{0\rho'}_{V',P}$, the moduli space in Proposition~\ref{intersect mdli}.\\

Below are two lemmas to embed every $M^{\rho'_i}_{V',P,ij}$ in $M_{V',P,0}$ for intersection purpose.

The base case is for $(\rho,H)$ in the case of Theorem~\ref{thm3}. With the same setting as in Theorem~\ref{thm3}. Let the morphism $M^{\rho}_{V,H}\xrightarrow[]{\iota} M_{V,H}$ be given by the canonical universal family representative over $M^{\rho}_{V,H}$ (see Remark~\ref{uni fml}). The morphism $\iota$ can be given explicitly by tracking the construction of both moduli spaces in Lemma~\ref{(M)}.
\begin{eg}\label{prt}
Here is an example that is a prototype for the morphism $M^{\rho}_{V,H} \rightarrow\mathbb{M}^{\rho}_{V,H} $ in the diagram of Lemma~\ref{(M)} below. The subring $k[X^\textsf{p}]$ of $k[X]$ is also a polynomial ring. The inclusion $k[X^\textsf{p}]\subset k[X]$ induces a bijection between closed points in $Spec(k[X])$ and those in $Spec(k[X^\textsf{p}])$, given explicitly by $(X-\lambda)\leftrightarrow(X^\textsf{p}-\lambda^\textsf{p})$.
\end{eg}

\begin{lem}\label{(M)}
There is a closed subscheme $\mathbb{M}^{\rho}_{V,H}$ of $M_{V,H}$ which $\iota$ factors through and whose closed points are in bijection with those of $M^{\rho}_{V,H}$ under $\iota$.
\[
\xymatrix@R=0.5cm{
  M^{\rho}_{V,H} \ar[dd]_{\iota} \ar[dr]^{}             \\
                & \mathbb{M}^{\rho}_{V,H} \ar@{^{(}->}[dl]_{}         \\
             M_{V,H}      }
\]
\end{lem}
\begin{proof}
Theorem~\ref{thm3}, Remark~\ref{dcps} and the base step for induction in the proof of Theorem 1.2 in [H80] are the references for this proof. Every fact used here can be found in one of the three places.

The explicit expression of $\iota$ on each $n$-th piece of $M^{\rho}_{V,H}$ (see Remark~\ref{clpt} c) will be given, using which the statements in the Lemma can be shown.

Denote by $M^{\rho}_{V,H,n}$ the $n$-th piece of $M^{\rho}_{V,H}$. The affine space $M^{\rho}_{V,H,n}$ can be identified with $Spec(k[K_n\spcheck-K_{n-1}\spcheck])$, where $K_n$, containing $K_{n-1}$, is the basis chosen for the $k$-vector space $(KerD)_n=KerD\cap H^0(V,q^n Div_{V})$ in the proof of Theorem~\ref{thm3}. Denote by $K_n\spcheck$ the set of the dual's of vectors in $K_n$. Write out elements in $K_n-K_{n-1}$ as $\{\verb"k"_i,1\leq i\leq d_K\}$. Then $Spec(k[K_n\spcheck-K_{n-1}\spcheck])=
Spec(k[\verb"k"_1\spcheck,...,
\verb"k"_{d_K}\spcheck])$.

Similarly denote by $M_{V,H,n}$ the $n$-th piece of $M_{V,H}$. The affine space $M_{V,H,n}$ can be identified with $Spec(k[L_n\spcheck-L_{n-1}\spcheck])$, where $L_n$, containing $L_{n-1}$, is the basis chosen for $H^0(V,q^nDiv_V)/{k^+}$. Denote by $L_n\spcheck$ the set of the dual's of vectors in $L_n$. Write out elements in $L_n-L_{n-1}$ as $\{l_j,1\leq j\leq d_L\}$. The way to choose $L_n$ is described in the base step for induction in the proof of Theorem 1.2 in [H80], analogous to the way to choose $K_n$. Only need to change the symbol $U$ there to $V$, $B_n$ there to $L_n$, and $A_n=H^0(U,q^nDiv_U)$ there to $B_n=H^0(V,q^nDiv_V)$. Recall that $U=Spec(A)$ and $V=Spec(B)$. $A_n$ and $B_n$ denote $k$-subspaces of $A$ and $B$ respectively.

Denote by $\iota_{n}$ the restriction of $\iota$ on $M^{\rho}_{V,H,n}$. The morphism $\iota_{n}$ maps every closed point in $M^{\rho}_{V,H,n}$ to the closed point in $M_{V,H,n}$ that represents the same pointed $H$-cover as it. Denote the $k$-algebra homomorphism that corresponds to $\iota_n$ by $\iota_n^*:k[L_n\spcheck-L_{n-1}\spcheck]\rightarrow k[K_n\spcheck-K_{n-1}\spcheck]$. It turns out that $\iota_n^*$ has the form: $l_j\spcheck \mapsto \Sigma_i\Sigma_{t\in X_j}(\lambda_{it}\verb"k"_i\spcheck)^{q_{tj}}$, where $X_j$ some finite set, $\lambda_{it}\in k$ and $q_{tj}$ is some $\mathsf{p}$-power.

The form of $\iota_n^*$ is obtained as follows. All pointed $H$-covers of $(V,v_g)$ can be given by elements in $B$ using Artin-Schreier equations $z^q-z=b$ with $b\in B$. Elements in the $k$-linear span of $L_n-L_{n-1}$ give bijectively all the pointed $H$-covers of $(V,v_g)$ that can be given by $z^q-z=b$ with $b\in B_n$. Every element $\sum_i \lambda_i\verb"k"_i$ in the $k$-linear span of $K_n -K_{n-1}$ is in $B_n$. Hence the pointed $H$-cover of $(V,v_g)$ given by $z^q-z=\sum_i \lambda_i\verb"k"_i$ is isomorphic to the pointed $H$-cover of $(V,v_g)$ given by $z^q-z=\sum_j \lambda'_jl_j$, for some unique $\sum_j \lambda'_jl_j$ in the $k$-linear span of $L_n-L_{n-1}$. The correspondence $\sum_i \lambda_i\verb"k"_i\leftrightarrow\sum_j \lambda'_jl_j$ is what is used to get the form of $\iota_n^*$: A closed point in $Spec(k[K_n\spcheck-K_{n-1}\spcheck])$ has the form $(\verb"k"_1\spcheck-\lambda_1,...,\verb"k"_{d_K}\spcheck-\lambda_{d_K})$. The maximal ideal represents the pointed $H$-cover of $(V,v_g)$ given by $z^q-z=\sum_i \lambda_i\verb"k"_i$, pointed anywhere above $v_g$. There is a unique $k$-algebra homomorphism $k[L_n\spcheck-L_{n-1}\spcheck]\rightarrow k[K_n\spcheck-K_{n-1}\spcheck]$ such that the inverse image of $(\verb"k"_1\spcheck-\lambda_1,...,\verb"k"_{d_K}\spcheck-\lambda_{d_K})$ is $(l_1\spcheck-\lambda'_1,...,l_{d_L}\spcheck-\lambda'_{d_L})$, which represents the pointed $H$-cover of $(V,v_g)$ given by $z^q-z=\sum_j \lambda'_jl_j$, for every closed point $(\verb"k"_1\spcheck-\lambda_1,...,\verb"k"_{d_K}\spcheck-\lambda_{d_K})$ in $Spec(k[K_n\spcheck-K_{n-1}\spcheck])$. Hence the homomorphism is $\iota_n^*$, by the definition of $\iota_n^*$. It is left as an exercise to the reader to write out the precise formula of the homomorphism, which has the form given above.

Let $\mathbb{M}^{\rho}_{V,H,n}=Spec(Im\iota_n^*)$, which is a closed subscheme of $M_{V,H,n}$. After simplification by elimination $Im\iota_n^*$ turns out a polynomial ring $k[\{\verb"k"'_i,1\leq i\leq d_K\}]$, where $\verb"k"'_i$ is a sum $\Sigma_{i<t\leq d_K} P_{it}(\verb"k"_t\spcheck)+
(\verb"k"_t\spcheck)^{n_{ii}}$ with $P_{it}$ a polynomial and $n_{ii}$ a $\mathsf{p}$-power. Moreover for every $i$ the polynomial ring $Im\iota_n^*$ contains a $(\verb"k"_i\spcheck)^{q_i}$ with $q_i$ a $\mathsf{p}$-power. Similar to Example~\ref{prt}, $\iota_n$ gives a bijection between the closed points of $M^{\rho}_{V,H,n}$ and those of $\mathbb{M}^{\rho}_{V,H,n}$.

The $\iota_n$'s for different $n$'s are compatible.
\end{proof}

Here are some necessary settings to prove the 2nd lemma for embedding $M^{\rho'_i}_{V',P,ij}$ in $M_{V',P,0}$.\\

With the same setting as in Theorem~\ref{thm2}. The ind scheme $M^{\rho,\bullet}_{V,P}$ consists of finitely many connected components \{$M^{\rho}_{V,P,j}$\}. For every $j$, the universal family (see Remark~\ref{unifmlpair} for more precise terminology) over $M^{\rho}_{V,P,j}$ determines a morphism $M^{\rho}_{V,P,j}\xrightarrow{\iota} M_{V,P}$, since $M_{V,P}$ is the fine moduli space for $F_{V,P}$.

If $(\rho_t,H_t)_t$ is a decomposition of $(\rho,P)$ (see Remark~\ref{dcps}), then $M^{\rho,\bullet}_{V,P}=\Pi_t M^{\rho_t,\bullet}_{V,H_t}$ and $M_{V,P}=\Pi_t M_{V,H_t}$. Hence $M^{\rho}_{V,P,j}$ has the form $\Pi_t M^{\rho_t}_{V,H_t}$ for every $j$.  The morphism $\iota$ can be given componentwise for each $t$.

\begin{lem}\label{cmap}
With the notations above, the morphism $M^{\rho}_{V,P,j}=\Pi_t M^{\rho_t}_{V,H_t}\xrightarrow{\iota}\Pi_t M_{V,H_t} $ is given by $\Pi_t \iota_t$, where $M^{\rho_t}_{V,H_t}\xrightarrow{\iota_t} M_{V,H_t}$ is the morphism given in Lemma~\ref{(M)}.
\end{lem}
\begin{proof}
Theorem~\ref{thm3}, Remark~\ref{dcps} and the base step for induction in the proof of Theorem 1.2 in [H80] are the references for this proof. Every fact used here can be found in one of the three places.

Induct on $|P|$.

(\ref{cmap}.1) The base case is done in Lemma~\ref{(M)}. Moreover for every $t$ since $M_{V,H_t}$ is the fine moduli space for $F_{V,H_t}$, the canonical universal family representative over $M_{V,H_t}$ given in 1.9 Rmk of [H80] pulled back to $M^{\rho_t}_{V,H_t}$ via $\iota_t$, differs from the canonical universal family representative over $M^{\rho_t}_{V,H_t}$ given in Remark~\ref{uni fml} by some element in $H^1(M^{\rho_t}_{V,H_t},H)$, by tracking definitions. Lemma~\ref{(M)} shows that $\iota_t$ gives a bijection on closed points of $M^{\rho_t}_{V,H_t}$ and $\mathbb{M}^{\rho_t}_{V,H_t}$. Using this fact and the same kind of argument in Lemma 4.25 of [TY17], a universal family representative over $M_{V,H_t}$ can be chosen such that it pulls back to the canonical universal family representative over $M^{\rho_t}_{V,H_t}$.

Below is the inductive step.

In Theorem~\ref{thm2} an $H$ inside the center $C(P)$ of $P$ is taken, and then the inductive process is carried out, which gives $M_{V,P}^{\rho\bullet}$ as $M^{\bar \rho\bullet}_{V,\bar P}\times M^{\rho_0\bullet}_{V,H}$. The notation $\Pi_t M^{\rho_t}_{V,H_t}$ means that $(H,\rho_0)$ is denoted by $(H_1,\rho_1)$ here and the inductive step there is carried out for some finite steps until the induction ends (see Remark~\ref{dcps}). Thus a connected component of $M_{V,\bar P}^{\bar\rho\bullet}$ can be denoted by $\Pi_{t\geq2} M^{\rho_t}_{V,H_t}$ and $M_{V,\bar P}$ by $\Pi_{t\geq2} M_{V,H_t}$.

(\ref{cmap}.2) By inductive hypothesis, the universal family over $\Pi_{t\geq2} M^{\rho_t}_{V,H_t}$ determines a morphism $\Pi_{t\geq2} M^{\rho_t}_{V,H_t}\xrightarrow{\Pi_{t\geq2} \iota_t}\Pi_{t\geq2} M_{V,H_t}$ with $\iota_t$ given in Lemma~\ref{(M)}; and there is a universal family representative over $\Pi_{t\geq2} M_{V,H_t}$ such that it pulls back to a $\bar\rho$-liftable universal family representative over $\Pi_{t\geq2} M^{\rho_t}_{V,H_t}$. Then lift the $\bar\rho$-liftable universal family representative over $\Pi_{t\geq2} M^{\rho_t}_{V,H_t}$ as in the proof of Theorem~\ref{thm2}, and pick any lift of the universal family representative over $\Pi_{t\geq2} M_{V,H_t}$. Again using the same kind of argument in Lemma 4.25 of [TY17], the latter lift can be modified such that its pullback to $\Pi_{t\geq2} M^{\rho_t}_{V,H_t}$ is the previous lift. (Strictly speaking, Definition~\ref{unifmlpair} needs to be used and pairs should be dealt with, which however will make the proof unnecessarily longer.)

A $\rho$-liftable universal family representative over $\Pi_t M^{\rho_t}_{V,H_t}$ can be given by the product of a $\rho_1$-liftable universal family representative over $M^{\rho_1}_{V,H_1}$ and a $ \rho$-liftable lift of a $\bar\rho$-liftable universal family representative over $\Pi_{t\geq2} M^{\rho_t}_{V,H_t}$ (see Remark~\ref{unifml2}). A universal family representative over $\Pi_t M_{V,H_t}$ is a similar product. Paragraphs (\ref{cmap}.1) and (\ref{cmap}.2) together show that the pullback of some universal family representative over $\Pi_t M_{V,H_t}$ via $\Pi_t M^{\rho_t}_{V,H_t}\xrightarrow{\Pi_t \iota_t}\Pi_t M_{V,H_t}$ is a $\rho$-liftable universal family representative over $\Pi_t M^{\rho_t}_{V,H_t}$. By the definition of the fine moduli space $M_{V,P}$,this means that $\Pi_t \iota_t$ is the morphism $\iota$ determined by the universal family over $M^{\rho}_{V,P,j}=\Pi_t M^{\rho_t}_{V,H_t}$.
\end{proof}

\begin{rem}\label{clsub}
By Lemma~\ref{(M)} and Lemma~\ref{cmap}, $\Pi_t \iota_t $ factors through $\Pi_t \mathbb{M}^{\rho_t}_{V,H_t}$, a closed subscheme of $M_{V,P}$.
\[
\xymatrix@R=0.5cm{
  \Pi_t M^{\rho_t}_{V,H_t} \ar[dd]_{\iota} \ar[dr]^{}             \\
                & \Pi_t \mathbb{M}^{\rho_t}_{V,H_t} \ar@{^{(}->}[dl]_{}         \\
             \Pi_t M_{V,H_t}      }
\]
\end{rem}

\begin{dfn}\label{attach}
If a $M_{V,P}=\Pi_t M_{V,H_t}$, a $M^{\rho}_{V,P,j}=\Pi_t M^{\rho_t}_{V,H_t}$ and their respective universal family representatives are constructed together, using the inductive process in the proof of Lemma~\ref{cmap}. Then the $M_{V,P}$ is called \textit{attached to} the $M^{\rho}_{V,P,j}$.
\end{dfn}
\begin{dfn}\label{M}
With the preparation of the two lemmas above, the moduli space $M^{0\rho'}_{V',P}$ for Proposition~\ref{intersect mdli} can be defined. There is a closed subscheme image (like $\Pi_t \mathbb{M}^{\rho_t}_{V,H_t}$ in Remark~\ref{clsub}) of $M^{\rho'_i}_{V',P,ij}$ in the $M_{V',P}$ attached (see Definition~\ref{attach}) to $M^{\rho'_i}_{V',P,ij}$. Using the isomorphism (see Remark~\ref{M's}) from the $M_{V',P}$ attached to $M^{\rho'_i}_{V',P,ij}$, to the fixed $M_{V',P,0}$, the closed subscheme image has its isomorphic image in $M_{V',P,0}$,  which is denoted by $\mathbb{M}^{\rho'_i}_{V',P,ij}$. Denote the morphism $M^{\rho'_i}_{V',P,ij}\rightarrow \mathbb{M}^{\rho'_i}_{V',P,ij}$ by $\iota_{\rho'_i,ij}$. Let $M^{\rho'}_{V',P}=\coprod_{(M^{\rho'_i}_{V',P,ij})_i}\bigcap_i \mathbb{M}^{\rho'_i}_{V',P,ij}$ (see between Remark~\ref{M's} and Example~\ref{prt} for the tuple $(M^{\rho'_i}_{V',P,ij})_i$). Let $M^0$ be the dense open subset of $M_{V',P,0}$ which parameterizes all connected pointed $P$-covers of $(V',v'_g)$ ([H80], Theorem 1.12). Let $M^{0\rho'}_{V',P}=\coprod_{(M^{\rho'_i}_{V',P,ij})_i}(
\bigcap_i \mathbb{M}^{\rho'_i}_{V',P,ij}
\bigcap M^0)$.
\end{dfn}
\begin{rem}\label{everyi}
What does the space $M^{0\rho'}_{V',P}$ parameterize?

Every closed point in $M^{0\rho'}_{V',P}$ represents a connected pointed $P$-cover $(W,w_g)\rightarrow (V',v'_g)$ corresponding to some homomorphism $\pi_1(V', v'_g)\xrightarrow{\phi}P$ that is $\rho'_i$-liftable for every $i$. In fact, for every $i$ the closed point gives a $\rho'_i$-liftable pair $(\phi,p_i)$ for some $p_i\in P$, by the definition of $M^{0\rho'}_{V',P}$ and the fact that every point in $M^{\rho\bullet}_{V,P}$ represents a $\rho$-liftable pair as shown in Theorem~\ref{thm2}.

See between Remark~\ref{M's} and Example~\ref{prt} for $c'_i$ and $\theta'$. The cover $(V',v'_g)\rightarrow(V'_{i'},\overline{v'_{gi}})$ corresponds to some homomorphism $\pi_1(V'_{i'},\overline{v'_{gi}})\xrightarrow{\theta'_i}\langle p_i'\rangle$ that maps $c'_i$ to $p_i'$. There is a similar diagram as in (\ref{glb mdli}.1) with $\widehat{\phi}_i(c'_i)=(p_i,p_i')$:
\[
  \xymatrix @R=.5cm @C=1.3cm {
  \pi_1(V', v'_g)\ar@{^{(}->}[d]_{}\ar[r]^{\phi} & P\ar@{^{(}->}[d]\\
  \pi_1(V'_{i'},\overline{v'_{gi}})\ar[r]^{\widehat{\phi}_i} & P\rtimes_{\rho_i'}\langle p_i'\rangle   \ar@{->>}[r]^{Q_P}& \langle p_i'\rangle
 ,}
 \]
where the composition of the bottom two arrows $\pi_1(V'_{i'},\overline{v'_{gi}})\rightarrow\langle p_i'\rangle$ is $\theta'_i$.

Hence the cover $(W,w_g)\rightarrow
(V'_{i'},\overline{v'_{gi}})$ is a pointed $P\rtimes_{\rho_i'}\langle p_i'\rangle$-cover. Denote by $\gamma_i$ the element in the Galois group of $W/V'_{i'}$ corresponding to $(1,p_i')$ with $1$ the identity of $P$. Denote by $\gamma_{p'_i}$ the automorphism in the Galois group of $V'\rightarrow U$ that corresponds to $p'_i\in P'$. Denote by $\gamma_p$ the element in the Galois group of $W/V'$ that corresponds to $p\in P$. $\gamma_i$ lies over $\gamma_{p'_i}$, and satisfies $ord(\gamma_i)=ord(p_i')$ and $\gamma_i\gamma_p {\gamma_i}^{-1}=\gamma_{\rho'(p_i')(p)}$ for every $p\in P$. These three conditions are called condition $(***i)$.

So a closed point in $M^{0\rho'}_{V',P}$ gives a pair $((W,w_g)\rightarrow (V',v'_g),\{\gamma_i\})$ with the first entry a connected pointed $P$-cover of $(V',v'_g)$ and the second entry a subset of the Galois group of $W/U$ with cardinality $|P'|$, the $i$-th element of which satisfies condition $(***i)$. The set of such pairs is denoted by $Gal_{V'}/\rho'$.

Reading backwards the discussion above, every pair in $Gal_{V'}/\rho'$ has a unique closed point in $M^{0\rho'}_{V',P}$ which represents the pair. Hence there is a canonical bijection between closed points in $M^{0\rho'}_{V',P}$ and $Gal_{V'}/\rho'$.
\end{rem}

\begin{rem}\label{gal gp}
There is a finite partition of closed points in $M^{0\rho'}_{V',P}$ by covers' Galois groups over $U$.

With the same notations as in the previous remark, the cover $W/U$ is Galois by a group order counting argument.

Denote $Gal(W/V')$ by $\Gamma_p$, $Gal(W/U)$ by $\Gamma$, and $Gal(V'/U)$ by $\Gamma_{p'} $. The isomorphism $\Gamma_p\simeq P$ is already given since $(W,w_g)\rightarrow (V',v'_g)$ is a $P$-cover. Similarly for $\Gamma_{p'}\simeq P'$.

Fix a subgroup $\widehat{\Gamma_{p'}}$ in $\Gamma$ which maps isomorphically to $\Gamma_{p'}$ under the canonical quotient map $\Gamma\twoheadrightarrow\Gamma_{p'}$. The existence of such a subgroup is given by Schur-Zassenhaus since $(\mathsf{p},|P'|)=1$. Then $\Gamma$ is canonically isomorphic to $\Gamma_p\rtimes\widehat{\Gamma_{p'}}$, an inner semiproduct.

The isomorphism $\Gamma_{p'}\simeq P'$ induces an isomorphism $\widehat{\Gamma_{p'}}\simeq P'$. Substituting $\Gamma_p$ by $P$ and $\widehat{\Gamma_{p'}}$ by $P'$ in $\Gamma_p\rtimes\widehat{\Gamma_{p'}}$, an induced semiproduct $P\rtimes_{\rho''}P'$ and an induced isomorphism $\Gamma\approx P\rtimes_{\rho''}P'$ are got. The diagram is commutative:
\[\xymatrix{
  1  \ar[r]^{} & \Gamma_p \ar[d]_{\simeq} \ar[r]^{} & \Gamma \ar[d]_{\approx} \ar[r]^{} & \Gamma_{p'} \ar[d]_{\simeq} \ar[r]^{} & 1  \\
  1 \ar[r]^{} & P \ar[r]^{} & P\rtimes_{\rho''}P' \ar[r]^{} & P' \ar[r]^{} & 1   .}\eqno{(\ref{gal gp}.1)}
  \]
For every $p'\in P'$, the action of $\rho'(p')$ on $P$ differs from that of $\rho''(p')$ by the conjugation of some element $p_{p'}\in P$, since $\gamma_i$ and its counterpart in $\widehat{\Gamma_{p'}}$ differ by some element in $\Gamma_p$.

When $P$ is abelian, the two groups $P\rtimes_{\rho'}P'$ and $P\rtimes_{\rho''}P'$ are the same. The Galois group over $U$ for any element in $Gal_{V'/\rho'}$ is $P\rtimes_{\rho'} P'$. If $P$ is not abelian, the two groups $P\rtimes_{\rho'}P'$ and $P\rtimes_{\rho''}P'$ may not be the same. The Galois group can not be nailed down.

The action $\rho''$ that arises above motivates the definition of a finite set consisting of certain semidirect products. Define $Gp_{\rho'}$ as the finite set \{$P\rtimes_{\rho''_s}P'\,|\,$ for every $p_i'$ there exists a $(p_i,p_i')$ in $P\rtimes_{\rho''_s}P'$ such that ord$(p_i,p_i')$=ord($p_i'$) and $(p_i,p_i')p(p_i,p_i')^{-1}=\rho'(p_i')(p)$\}, where $\rho''_s:P'\rightarrow Aut(P)$ is an action of $P'$ on $P$.

By the end of Remark~\ref{everyi}, the closed points of $M^{0\rho'}_{V',P}$ are in canonical bijection with $Gal_{V'}/{\rho'}$. Every pair in $Gal_{V'}/{\rho'}$ gives a pointed $P\rtimes_{\rho''_s}P'$-cover (a similar diagram to diagram (\ref{glb mdli}.1)):
\[
  \xymatrix @R=.5cm @C=1.3cm {
  \pi_1(V', v'_g)\ar@{^{(}->}[d]_{}\ar[r]^{\phi} & P\ar@{^{(}->}[d]\\
  \pi_1(U, u_g)\ar[r]^{\widehat{\phi}} &  P\rtimes_{\rho''_s}P'  \ar@{->>}[r]^{Q_P}& P'
 ,} \eqno{(\ref{gal gp}.2)}
 \]
for some $\rho''_s$ using the process given above diagram (\ref{gal gp}.1), where the composition of the bottom two arrows is $\theta'$. The group $P\rtimes_{\rho''_s}P'$ is said to \textit{belong to the pair or belong to the closed point corresponding to the pair}. A different $P\rtimes_{\rho''_{s_1}}P'$ can belong to the same pair, if a different section $\widehat{\Gamma_{p'}}$ is chosen in the process. If two $P\rtimes_{\rho''_s}P'$ and $P\rtimes_{\rho''_{s_1}}P'$ belong to the same pair, then a similar diagram to (\ref{gal gp}.2)
\[
  \xymatrix @R=.5cm @C=1.3cm {
  \pi_1(V', v'_g)\ar@{^{(}->}[d]_{}\ar[r]^{\phi} & P\ar@{^{(}->}[d]\\
  \pi_1(U, u_g)\ar[r]^{\widehat{\phi}_1} &  P\rtimes_{\rho''_{s_1}}P'  \ar@{->>}[r]^{Q_P}& P'
 }
 \]
is also commutative. It together with diagram (\ref{gal gp}.2) gives a commutative diagram:
\[\xymatrix{
  1  \ar[r]^{} & P \ar[d]_{=} \ar[r]^{} & P\rtimes_{\rho''_{s_1}}P' \ar[d]_{\simeq} \ar[r]^{} & P' \ar[d]_{=} \ar[r]^{} & 1  \\
  1 \ar[r]^{} & P \ar[r]^{} & P\rtimes_{\rho''_s}P' \ar[r]^{} & P' \ar[r]^{} & 1   .}
  \]
Hence $P\rtimes_{\rho''_s}P'$ and $P\rtimes_{\rho''_{s_1}}P'$ are isomorphic extensions. Pick a representative (pick $P\rtimes_{\rho'}P'$ in its class) in each isomorphism class of extensions and denote the subset obtained in this way of $Gp_{\rho'}$ by $\bar{Gp}_{\rho'}$. The set $Gal_{V'}/{\rho'}$ has a finite partition by elements in $\bar{Gp}_{\rho'}=
\{P\rtimes_{\rho''_t}P'\}$ by discussion above.

For any $P\rtimes_{\rho''_s}P'\in Gp_{\rho'}$ and any pointed $P\rtimes_{\rho''_s}P'$-cover $(W,w_g)\rightarrow(U,u_g)$ corresponding to some $\pi_1(U,u_g) \xrightarrow{ \hat{ \hat{\phi} } }   P\rtimes_{\rho''_s}P'$, a pointed $P\rtimes_{\rho''_{si}}\langle p_i'\rangle$-cover $(W,w_g)\rightarrow(V'_{i'},
\overline{v'_{gi}})$ can be got for every $i$:
\[\xymatrix{
  \pi_1(V',v'_g)\ar@{^{(}->}[d]_{} \ar[r]^{} & P \ar@{^{(}->}[d]^{}\ar[r]^{=} & P\ar@{^{(}->}[d]^{} \\
  \pi_1(V'_{i'},\overline{v'_{gi}}) \ar@{->}[d]_{} \ar[r]^{} & P\rtimes_{\rho''_{si}}\langle p_i'\rangle \ar[r]^{\simeq}\ar@{^{(}->}[d]^{} & P\rtimes_{\rho'_{i}}\langle p_i'\rangle \\
  \pi_1(U,u_g) \ar[r]^{\hat{\hat{\phi}}} & P\rtimes_{\rho''_s}P'   ,}\eqno{(\ref{gal gp}.3)}\]
where $\rho''_{si}$ is the restriction of $\rho''_s$ on $\langle p_i'\rangle$.
The pointed $P\rtimes_{\rho''_{si}}\langle p_i'\rangle$-cover $(W,w_g)\rightarrow(V'_{i'},\overline{v'_{gi}})$ is also a pointed $P\rtimes_{\rho'_{i}}\langle p_i'\rangle$-cover, as shown in the commutative diagram (\ref{gal gp}.3), where the group isomorphism $P\rtimes_{\rho''_{si}}\langle p_i'\rangle \rightarrow P\rtimes_{\rho'_{i}}\langle p_i'\rangle $ sends $(p_i,p_i')$ to $(1,p_i')$ and every $p\in P$ to $p$. Then Remark~\ref{everyi} shows that $(W,w_g)\rightarrow(U,u_g)$ gives a pair in $Gal_{V'}/{\rho'}$ corresponding to some closed point in $M^{0\rho'}_{V',P}$, which can be used to discover the original $(W,w_g)\rightarrow(U,u_g)$ using diagram (\ref{gal gp}.3). If a closed point in $M^{0\rho'}_{V',P}$ is used to discover, using diagram (\ref{gal gp}.3), a pointed $P\rtimes_{\rho''_s}P'$-cover of $(U,u_g)$ for some $P\rtimes_{\rho''_s}P'$ in $Gp_{\rho'}$, there are several possibilities for $\rho''_s$.
\end{rem}

To define the functor $F^{er,Gal/\rho'}_{V',P}$ in Proposition~\ref{intersect mdli}, several new definitions are needed. The inclusion of polynomial rings $Im\iota_{tn}^*\subset k[\{\verb"k"_i,1\leq i\leq d_K\}]$ in the proof of Lemma~\ref{(M)} and the morphism $\Pi_t M^{\rho_t}_{V,H_t} \rightarrow\Pi_t \mathbb{M}^{\rho_t}_{V,H_t}$ in Remark~\ref{clsub} motivate the first two definitions given below respectively.
\begin{dfn}\label{type-R}
Let $k[X_1,...,X_d]$ be a polynomial ring. Suppose $P'_0$ is a subring that can be written as, for every permutation $\{s(i)\}$ of $\{1,...,d\}$, a polynomial ring $k[X'_1,...,X'_d]$ with each $X'_i$ a sum $\Sigma_{i< t\leq d} P_{it}(X_{s(t)})+(X_{s(i)})^{n_{ii}}$ , $P_{it}$ a polynomial and $n_{ii}$ a $\mathsf{p}$-power. Let $\mathcal{P}'$ be a polynomial ring with an injective $k$-algebra homomorphism $f:\mathcal{P}'\hookrightarrow k[X_1,...,X_d]$. If $f$ gives an isomorphism between $\mathcal{P}'$ and $P'_0$, then $f:\mathcal{P}'\hookrightarrow k[X_1,...,X_d]$ is an \textit{$R$-extension}.

Let $\{\mathcal{P}_i\hookleftarrow \mathcal{P}_i'\}$ be a collection of finitely many $R$-extensions, with possibly different $\mathcal{P}_i$'s and $\mathcal{P}'_i$'s. Tensoring over $k$ gives a morphism $Spec(\otimes_i\mathcal{P}_i) \rightarrow Spec(\otimes_i\mathcal{P}'_i) $. For any $(S',s'_0)\in \mathcal{S}$ and $(S',s'_0)\xrightarrow[]{f} (Spec(\otimes_i\mathcal{P}'_i),x'_0)$ a morphism in $\mathcal{S}$, the pullback $(S,s_0)\rightarrow (S',s'_0)$ of $(Spec(\otimes_i\mathcal{P}_i),x_0)\rightarrow (Spec(\otimes_i\mathcal{P}'_i),x'_0)$, for some $x_0$ mapping to $x'_0$, is called \textit{a morphism of type $R$}:
 \[
  \xymatrix @R=.5cm @C=1.8cm {
  S\ar@{->}[d]_{}\ar[r]^{f'} & Spec(\otimes_i\mathcal{P}_i)\ar@{->}[d]\\
  S'\ar[r]^{f} &  Spec(\otimes_i\mathcal{P}'_i)
 .}
 \]

A morphism $(T,t_0)\rightarrow(S,s_0)$ in $\mathcal{S}$ is of \textit{type $ER$}, if it can be decomposed into a finite sequence of finite etale covers and morphisms of type $R$.
\end{dfn}
\begin{rem}\label{MM}
The morphism $M^{\rho'_i}_{V',P,ij}\xrightarrow{\iota_{\rho'_i,ij}} \mathbb{M}^{\rho'_i}_{V',P,ij}$ in Definition~\ref{M} is an example of the right column in the square diagram in Definition~\ref{type-R}.
\end{rem}

Below is the definition for the functor to present in Proposition~\ref{intersect mdli}, which is motivated by the discussion in Remark~\ref{everyi}. Let $(S,s_0)\in\mathcal{S}$ and let $Gal_{(S,s_0)}$ be the set of $T$-parameterized $P$-covers $(\widetilde W,\tilde{w}_g)\rightarrow (T\times V',(t_0,v'_g))$ of $V'$ pointed over $(t_{0},v'_g)$ for some $(T,t_{0})\rightarrow (S,s_0)$ of type $ER$ with connected fibers over the closed points of $T$, such that the composition $(\widetilde W,\tilde{w}_g)\rightarrow (T\times V',(t_0,v'_g))\rightarrow (T\times U,(t_0,u_g))$ is Galois. For a pointed $P$-cover $(\widetilde W,\tilde{w}_g)\rightarrow (T\times V',(t_0,v'_g))$, and an element $p\in P$, denote by $\tilde{\gamma}_p$ the automorphism in its Galois group that corresponds to $p$. Denote by $\tilde{\gamma}_{p'_i}$ the automorphism in the Galois group of $T\times V'\rightarrow T\times U$ that corresponds to $p'_i\in P'$. Let $Gal/\rho'_{(S,s_0)}$ be the set of pairs $((\widetilde W,\tilde{w}_g)\rightarrow (T\times V',(t_0,v'_g))$,\{$\tilde{\gamma}_i$\}), where $(\widetilde W,\tilde{w}_g)\rightarrow (T\times V',(t_0,v'_g))$ is in $Gal_{(S,s_0)}$ and $\tilde{\gamma}_i$ is in the Galois group of the cover $\widetilde W\rightarrow T\times U$ that lies over $\tilde{\gamma}_{p'_i}$ and satisfies $ord(\tilde{\gamma}_i)=ord(p_i')$ and $\tilde{\gamma}_i\tilde{\gamma}_p {\tilde{\gamma}_i}^{-1}=\tilde{\gamma}_{\rho'(p_i')(p)}$ for every $p\in P$. The set $Gal/\rho'_{(S,s_0)}$ is an $S$-parameterized version of $Gal_{V'/\rho'}$; by Remark~\ref{everyi}, the set $F^{er,Gal/\rho'}_{V',P}(Spec(k),s_0)$, with $s_0$ determined by $v'_g$ using diagram (\ref{secn&t}.1), is the set $Gal_{V'}/\rho'$. Two elements ($(\widetilde W_j,\tilde{w}_{gj})\rightarrow (T_j\times V',(t_{j0},v'_g))$,\{$\tilde\gamma_{ji}$\}) ($j=1,2$) in $Gal/\rho'_{(S,s_0)}$ are \textit{$ER$-equivalent} if there exists a morphism $(T_d,t_{d0})\rightarrow (S,s_0)$ of type $ER$, where $(T_d,t_{d0})$ also maps to $(T_j,t_{j0})$ ($j=1,2$) in the category $\mathcal{S}$, such that the two pointed $P$-covers, together with \{$\tilde\gamma_{1i}$\} and \{$\tilde\gamma_{2i}$\}, pulled back to $T_d$ become isomorphic. Let $F^{er,Gal/\rho'}_{V',P}$ be the functor: $\mathcal{S}\rightarrow$ (Sets); $(S,s_0)\mapsto$ \{$ER$-equivalence classes of $((\widetilde W,\tilde{w}_g)\rightarrow (T\times V',(t_0,v'_g))$,\{$\tilde\gamma_i$\}) $\in Gal/\rho'_{(S,s_0)}$\}.

Here is the last definition involved in the statement of Proposition~\ref{intersect mdli}. Two morphisms $T_j\xrightarrow[]{f_j} M^{0\rho'}_{V',P}$ ($j=$1 or 2, where $(T_j,t_{j0})\rightarrow (S,s_0)$ is of type $ER$) are \textit{$ER$-equivalent}, if there exists a morphism $(T_d,t_{d0})\rightarrow (S,s_0)$ of type $ER$ with $(T_d,t_{d0})$ also mapping to $(T_j,t_{j0})$ ($j=1,2$) in the category $\mathcal{S}$, such that the $f_j$'s pulled back to $T_d$ are the same.  Let $M^{er0\rho'}_{V',P}$ be the functor: $\mathcal{S}$ $\rightarrow$ (Sets); $(S,s_0)\mapsto$\{$ER$-equivalence classes of $T\rightarrow M^{0\rho'}_{V',P}$, where $(T,t_0)\rightarrow (S,s_0)$ runs over all morphisms to $S$ of type $ER$\}.
\begin{rem}\label{er}
The two $ER$-equivalences in the definitions of functors $F^{er,Gal/\rho'}_{V',P}$ and $M^{er0\rho'}_{V',P}$ arise naturally in the proof of Proposition~\ref{intersect mdli} based on the intersection idea.
\end{rem}
\begin{pro}\label{intersect mdli}
With the same notations as above, the ind scheme $M^{0\rho'}_{V',P}$ is the moduli space for $F^{er,Gal/\rho'}_{V',P}$ in the sense that there exists an isomorphism between functors $F^{er,Gal/\rho'}_{V',P}\simeq M^{er0\rho'}_{V',P}$.

Moreover, on each of the finitely many irreducible components of $M^{0\rho'}_{V',P}$, there is a unique $P\rtimes_{\rho''_t}P'$ in $\bar{Gp}_{\rho'}$ which belongs to (defined in Remark~\ref{gal gp}) all the closed points. Conversely, for every $P\rtimes_{\rho''_t}P'$ in $\bar{Gp}_{\rho'}$, there is an irreducible component, such that $P\rtimes_{\rho''_t}P'$ belongs to all the closed points of the component.
\end{pro}
\begin{proof}
Proof of the first statement:

Let $(S,s_0)\in \mathcal{S}$ and $((\widetilde W,\tilde w_g)\rightarrow (T\times V',(t_0,v'_g)),\{\gamma_i\})$ be a representative in an $ER$-equivalence class of $F^{er,Gal/\rho'}_{V',P}(S,s_0)$. Then $(\widetilde W,\tilde w_g)\rightarrow (T\times V'_{i'}, (t_0,\overline{v'_{gi}}))$ is Galois (see Remark~\ref{everyi}). Letting $\gamma_i$ correspond to $(1,p_i')\in P\rtimes_{\rho'_i}\langle p_i'\rangle$, $(\widetilde W,\tilde w_g)\rightarrow (T\times V'_{i'}, (t_0,\overline{v'_{gi}}))$ is a pointed $P\rtimes_{\rho'_i}\langle p_i'\rangle$-cover. By the definition of $M^{\rho'_i,\bullet}_{V',P}$, the pointed $P\rtimes_{\rho'_i}\langle p_i'\rangle$-cover corresponds to a morphism $T\xrightarrow[]{\mathfrak{c}_i} M^{\rho'_i,\bullet}_{V',P}$. Since $T$ is connected, the morphism $\mathfrak{c}_i$ lands in a connected component $M^{\rho'_i}_{V',P,ij}$ of $M^{\rho'_i,\bullet}_{V',P}$. Embedding $M^{\rho'_i}_{V',P,ij}$ in $M_{V',P,0}$ as in Remark~\ref{MM}, a morphism $T\xrightarrow[]{\hat{\mathfrak{c}}_i}M_{V',P,0}$ is got. The morphism $\hat{\mathfrak{c}}_i$ is the same as the morphism from $T$ to $M_{V',P,0}$ determined by the pointed $P$-cover $(\widetilde W,\tilde w_g)\rightarrow (T\times V',(t_0,v'_g))$, using that $M_{V',P,0}$ is the fine moduli space for  pointed families of $P$-covers of $(V',v'_g)$ ([H80], Theorem 1.2). The above discussion applies for every $i$. Hence a morphism $T\xrightarrow[]{\hat{\mathfrak{c}}}M^{0\rho'}_{V',P}$ is got, by the definition of $M^{0\rho'}_{V',P}$.

Conversely, given $T\xrightarrow[]{ \hat{\mathfrak{c}} }M^{0\rho'}_{V',P}$ for some $(T,t_{0})\rightarrow (S,s_0)$ of type $ER$, a morphism $T\xrightarrow[]{ \hat{\mathfrak{c}} }M_{V',P,0}$ is got by the definition of $M^{0\rho'}_{V',P}$ and the connectedness of $T$. For each $i$, there is a morphism $T\xrightarrow[]{\hat{\mathfrak{c}}_i}
\mathbb{M}^{\rho'_i}_{V',P,ij}$, for some $\mathbb{M}^{\rho'_i}_{V',P,ij}$ (see Definition~\ref{M} for $\mathbb{M}^{\rho'_i}_{V',P,ij}$) that $\hat{\mathfrak{c}}$ factors through
\[
\xymatrix{
  T \ar[dr]_{\hat{\mathfrak{c}}}\ar[r]^{\hat{\mathfrak{c}}_i}
                &  \mathbb{M}^{\rho'_i}_{V',P,ij}\ar@{^{(}->}[d]^{}  \\
                & M_{V',P,0}            .}
\]

After pointing $M^{\rho'_i}_{V',P,ij}\xrightarrow
{\iota_{\rho'_i,ij}}
\mathbb{M}^{\rho'_i}_{V',P,ij}$ (see Definition~\ref{M}) properly, the pullback morphism $(T_i,t_{i0})\rightarrow (T,t_{0})$ is of type $R$, by Remark~\ref{MM} and the lower square of the diagram, in which both squares are pullbacks:
\[
  \xymatrix @R=.5cm @C=1.3cm {
  T_{di}\ar[d]_{}\ar[r]^{} & M_{ij}\ar@{->}[d]\\
  T_i\ar[d]_{}\ar[r]^{\mathfrak{c}_i} & M^{\rho'_i}_{V',P,ij}\ar[d]\\
  T\ar[r]^{\hat{\mathfrak{c}}_i} &  \mathbb{M}^{\rho'_i}_{V',P,ij}
 .}
 \]
The (ind) (see Remark~\ref{n}) scheme $M_{ij}$ in the upper right corner is a finite etale cover of $M^{\rho'_i}_{V',P,ij}$, such that a fixed universal family representative over $M_{V',P,0}$ pulled back to $M_{ij}$, is the same as the pullback to $M_{ij}$ of a $\rho'_i$-liftable universal family representative over $M^{\rho'_i}_{V',P,ij}$ (see the discussion between Lemma~\ref{(M)} and Lemma~\ref{cmap}).

The relationship between $M_{ij}$ and $M^{\rho'_i}_{V',P,ij}$ is the same as that between $T$ and $S$ in the proof of Theorem~\ref{thm2}. Hence the pullback $(T_{di},t_{di0})\rightarrow (T_i,t_{i0})$ is a finite etale cover.

The two diagrams together imply that the fixed universal family representative over $M_{V',P,0}$, which is a pointed $P$-cover of $(M_{V',P,0}\times V',(\hat{\mathfrak{c}}(t_0),v'_g))$, pulled back to $T_{di}$ is $\rho'_i$-liftable. Let $(T_d,t_{d0})$ be the common pullback of the $(T_{di},t_{di0})$'s over $(T,t_0)$. The pullback to $T_d$ of the fixed universal family representative over $M_{V',P,0}$ is $\rho'_i$-liftable for every $i$. Consider $\rho'_i$-liftable pairs, rather than merely $\rho'_i$-liftable covers, at some places in the discussion above. Then a pair (see Remark~\ref{everyi}) $((\widetilde W,\tilde w_g)\rightarrow (T_d\times V',(t_{d0},v'_g)),\{\gamma_i\})$ is got, a pointed $P$-cover together with $|P'|$ elements in $Gal(\widetilde W/T_d\times U)$, whose $ER$-equivalence class is in $F^{er,Gal/\rho'}_{V',P}(S,s_0)$.

The two maps are well defined for equivalence classes and inverse to each other.

Proof of the statements after ``Moreover": Every component $\bigcap_i \mathbb{M}^{\rho'_i}_{V',P,ij}
\bigcap M^0$ of $M^{0\rho'}_{V',P}$ is a dense open of $\bigcap_i \mathbb{M}^{\rho'_i}_{V',P,ij}$, which itself is an affine closed subscheme of  $M_{V',P,0}$. Pick any irreducible component of $\bigcap_i \mathbb{M}^{\rho'_i}_{V',P,ij}
\bigcap M^0$ and a covering of it consisting of connected affine open subsets of finite type over $k$, which are all dense and intersect each other. Denote any of the affine open subsets by $\mathbb{M}$ and apply the first statement proven above to $\mathbb{M}$. The inclusion of $\mathbb{M}$ in $M^{0\rho'}_{V',P}$, with any base point $m_g$, gives a pointed $P$-cover of $(\mathbb{M}'\times V',(m'_g,v'_g))$ for some $(\mathbb{M}',m'_g)\rightarrow(\mathbb{M},m_g)$ of type-$ER$. The pointed $P$-cover satisfies a $\mathbb{M}'$-parameterized version of $(***i)$ for every $i$ and thus gives a pointed $P\rtimes_{\rho''_t}P'$-cover of $(\mathbb{M}'\times U,(m'_g,u_g))$ for some unique $P\rtimes_{\rho''_t}P'$ in $\bar{Gp}_{\rho'}$ (see Remark~\ref{everyi} and Remark~\ref{gal gp}). Hence for every closed point $m$ (need to use chemins for base point issues) in $\mathbb{M}$, $P\rtimes_{\rho''_t}P'$ belongs to the pair in $Gal_{V'/\rho'}$ that $m$ corresponds. Then Remark~\ref{everyi} and Remark~\ref{gal gp} suffice to give all the statements.
\end{proof}
Denote the maximal union of irreducible components of $M^{0\rho'}_{V',P}$, to all of whose closed points $P\rtimes_{\rho'}P'$ belongs (see the last paragraph in the proof of Proposition~\ref{intersect mdli}), by $M_{U,V',P\rtimes_{\rho'}P'}$. Denote by $M_{U,P\rtimes_{\rho'}P'}$ the disjoint union of $M_{U,V',P\rtimes_{\rho'}P'}$'s over all possible $(V',v'_g)$'s pointed connected $P'$-covers of $(U,u_g)$.

Define a functor $M^{er}_{U,P\rtimes_{\rho'}P'}$: $\mathcal{S}$ $\rightarrow$ (Sets); $(S,s_0)\mapsto$\{$ER$-equivalence classes of $T\rightarrow M_{U,P\rtimes_{\rho'}P'}$, where $(T,t_0)\rightarrow (S,s_0)$ runs over all morphisms to $S$ of type $ER$\}. Similarly to $M^{er0\rho'}_{V',P}$ and $M^{0\rho'}_{V',P}$ above.

Define a functor $F^{er}_{U,P\rtimes_{\rho'}P'}$: $\mathcal{S}$ $\rightarrow$ (Sets); $(S,s_0)\mapsto$\{$ER$-equivalence classes of pointed $P\rtimes_{\rho'}P'$-covers $(\widetilde{W},\tilde{w}_g)\rightarrow(T\times U,(t_0,u_g))$ whose fibers over closed points of $T$ are all connected, where $(T,t_0)\rightarrow (S,s_0)$ runs over all morphisms to $S$ of type $ER$\}.  The definition of $ER$-equivalence classes here is obvious (see definition of the functor $F^{er,Gal/\rho'}_{V',P}$).
\begin{cor}\label{uniirr}
The functor $F^{er}_{U,P\rtimes_{\rho'}P'}$ is isomorphic to the functor $M^{er}_{U,P\rtimes_{\rho'}P'}$, which shows that $M_{U,P\rtimes_{\rho'}P'}$ is a moduli space for $P\rtimes_{\rho'}P'$-covers of $(U,u_g)$.
\end{cor}
\begin{proof}
Directly from Proposition~\ref{intersect mdli} and its proof. See also proof of Theorem~\ref{thm1}.
\end{proof}

\section{Local vs. global moduli} \label{l vs g}
In Section~\ref{l vs g}, a fine moduli space (Proposition~\ref{ram mdli-3}) for cyclic-by-$\mathsf{p}$ covers of an affine curve at most tamely ramified over finitely many closed points, is constructed. The new type of fine moduli space is obtained by modifying the proof for the previous global fine moduli space constructed in Theorem~\ref{thm1} Section~\ref{glb mdli}, and is constructed in similar 3 steps. The new type of fine moduli space is the global side of a local-global principal Proposition~\ref{local-global moduli}. There is a different phenomenon for cyclic-by-$\mathsf{p}$ covers from that for $\mathsf{p}$-covers. In [H80] the similar local-global principal for $\mathsf{p}$-groups stated in Proposition 2.1 does not involve (tamely) ramified global covers; there the global covers are etale. The local-global principal Proposition~\ref{local-global moduli} has a version over a general field of characteristic $\mathsf{p}>0$, which is Main Theorem 1.4.1 in [K86].

A parameter space for local cyclic-by-$\mathsf{p}$ covers of $Spec(k((x)))$ is constructed in Proposition~\ref{local para space-3}, which is the local side of the local-global principal Proposition~\ref{local-global moduli}. The construction is also by modifying the one in Section~\ref{glb mdli} and has similar 3 steps.

Finally it is shown that a restriction morphism (a general case of the local-global principal Proposition~\ref{local-global moduli} with the isomorphism there replaced by a finite morphism now) is finite, which is from the new type of global moduli space to a product of the local parameter spaces (Proposition~\ref{fin etl-1}), an analogue to Proposition 2.7 in [H80]. It is proved by a similar argument.

As always, we follow notations and terminology defined in Section~\ref{secn&t}. For example $G$ represents a cyclic-by-$\mathsf{p}$ group.\\

Here are some necessary settings for the construction of the fine moduli space (Proposition~\ref{ram mdli-3}).

Let $T$ be a finite set of closed points on $U$ not including $u_g$ and $U^0=U-T$. Denote by \{$(V^0_l,v_l)$\} the set of all the finitely many connected pointed $\mathbb{Z}/n_l$-covers of $(U^0,u_g)$, where $n_l$ can be any factor of $n$. Let $(V_l,v_l)\rightarrow(U,u_g)$ be the extension of $(V^0_l,v_l)\rightarrow(U^0,u_g)$, obtained by putting back in some deleted closed points from the smooth completions of both curves.

\begin{dfn}\label{F1T}
Let $F^T_{U,G}$ be the functor: $\mathcal{S}_1\rightarrow$(Sets), $(S,s_0)\mapsto$\{equivalence classes of possibly ramified $G$-covers $\widetilde W\rightarrow S\times U$ pointed over $(s_0,u_g)$, where the restriction of $\widetilde W$ over $S\times U^0$ is a $G$-cover and $\widetilde W\rightarrow \widetilde W/P$ is finite etale\}. Let $S=Spec(k)$ with $s_0$ determined by $u_g$ using diagram (\ref{secn&t}.1). Then $F^T_{U,G}(S,s_0)$ is the set of possibly ramified pointed $G$-covers $ (W,w_g)\rightarrow(U,u_g)$ whose restriction $W^0$ over $U^0$ is a $G$-cover and $W\rightarrow W/P$ is finite etale.
\end{dfn}

A group homomorphism $\widetilde\phi^0:\pi_1(S\times V_l^0,(s_0,v_l))\rightarrow P$ \textit{factors through $V_l$} if $\widetilde\phi^0=(\pi_1(S\times V_l^0,(s_0,v_l))\xrightarrow{(Id_S\times (V_l^0\subset V_l) )_*}\pi_1(S\times V_l,(s_0,v_l))\xrightarrow{\widetilde{\phi}} P$) for some $\widetilde{\phi}$.

\begin{dfn}\label{F2T}
Let $F^{\rho_{n_l}\bullet/T}_{V_l^0,P}$ be the functor: $\mathcal{S}_1\rightarrow$(Sets), $(S,s_0)\mapsto$ \{([$\widetilde\phi^0$], $p$), where $(\widetilde\phi^0,p)$ is a $\rho_{n_l}$-liftable pair with $\widetilde\phi^0$ factoring through $V_l$.\}. Let $S=Spec(k)$ with $s_0$ determined by $v_l$ using diagram (\ref{secn&t}.1). Then $F^{\rho_{n_l}\bullet/T}_{V_l^0,P}
(S,s_0)$ is the set of $\rho_{n_l}$-liftable pairs of $(V_l^0,v_l)$ the first entries of which can all extend to $P$-covers of $V_l$.
\end{dfn}
Below is the first of the 3 steps in constructing the fine moduli space in Proposition~\ref{ram mdli-3}.
\begin{dfn}\label{newlftb}(Remark/Definition)

The cover $V_l\rightarrow U$ above is ramified at finitely many closed points on $U$. Hence the old definition of $\rho$-liftable can not apply here. New definition: A group homomorphism $\phi:\pi_1(V_l,v_l)\rightarrow H$ with $(\rho,H)$ in the case of Theorem~\ref{thm3} is \textit{$\rho$-liftable} if $\phi$ makes the diagram in Lemma~\ref{ff} commutative. Using the new definition of $\rho$-liftable in the definition of $F_{V,H}^{\rho}$, Theorem~\ref{thm3} still holds and $F_{V,H}^{\rho}$ is presented by the same fine moduli space $M_{V,H}^{\rho}$.

The new definition of $\rho$-liftable is used below for $V_l\rightarrow U$.
\end{dfn}
\begin{lem}\label{ram mdli-1}
With the notations above, suppose $(\rho_{n_l},P)$ is in the case of Theorem~\ref{thm3} and write $H$ for $P$ in this case. The functor $F^{\rho_{n_l}\bullet/T}_{V_l^0,H}$ for $\rho_{n_l}$-liftable pairs of $(V_l^0,v_l)$ the first entries of which can all extend to $H$-covers of $V_l$, has a fine moduli space, a disjoint union of finitely many ind affine spaces.
\end{lem}
\begin{proof}
It is to be shown that $F^{\rho_{n_l}\bullet/T}_{V_l^0,H}=\amalg_h F^{\rho_{n_l}}_{V_l,H}$, copies of $F^{\rho_{n_l}}_{V_l,H}$ indexed by $h$, where depending on $(\rho_{n_l},H)$ the $h$ runs over $H$ or it is just $1$ (see Corollary~\ref{disj uni over H}). Suppose a $\rho_l$-liftable pair $(\widetilde{\phi^0}:\pi_1(S\times V_l^0,(s_0,v_l))\rightarrow H,h_0)$ with $\widetilde{\phi^0}$ factoring through $V_l$ is a representative for an equivalence class ($[\widetilde{\phi^0}],h_0$) in $F^{\rho_{n_l}\bullet/T}_{V_l^0,H}(S,s_0)$. Since $(Id_S\times (V_l^0\subset V_l) )_*$ is surjective, $\widetilde{\phi}$ is also $\rho_{n_l}$-liftable (see Definition~\ref{newlftb} and see the definition of ``factoring through $V_l$" above for $\widetilde{\phi}$). The left to right map sends ($[\widetilde{\phi^0}],h_0$) to $[\widetilde{\phi}]$ in the copy of $F^{\rho_{n_l}}_{V_l,H}(S,s_0)$ indexed by $h_0$. The inverse map is obvious. By Theorem~\ref{thm3} (see also Definition~\ref{newlftb}), $F^{\rho_{n_l}}_{V_l,H}$ is represented by an ind affine space $M^{\rho_{n_l}}_{V_l,H}$.
\end{proof}

Here is the 2nd of the 3 steps in constructing the fine moduli space in Proposition~\ref{ram mdli-3}.
\begin{lem}\label{ram mdli-2}
The functor $F^{\rho_{n_l}\bullet/T}_{V_l^0,P}$ for $\rho_{n_l}$-liftable pairs of $(V_l^0,v_l)$ the first entries of which can all extend to $P$-covers of $V_l$, has a fine moduli space, a disjoint union of finitely many ind affine spaces.
\end{lem}
\begin{proof}
The proof is by simply replacing with their obvious counterparts symbols in and slightly modifying the proof of Theorem~\ref{thm2}.

Denote $V_l$ by $V$ and $\rho_{n_l}$ by $\rho$ in the proof. Replace $F^{\rho\bullet}_{V,P}$ with $F^{\rho\bullet/T}_{V^0,P}$, Corollary~\ref{disj uni over H} with Lemma~\ref{ram mdli-1}, $M^{\bar \rho\bullet}_{V,\bar P}$ with $M^{\bar \rho\bullet/T}_{V^0,\bar P}$, and $M^{\rho_0\bullet}_{V,H}$ with $M^{\rho_0\bullet/T}_{V^0,H}$.

In paragraph (\ref{thm2}.2) of the proof of Theorem~\ref{thm2}, since the first entry of some universal pair ($\widetilde{\mu}^0_{0}: \pi_1(\bar M\times V^0,(\bar{m},v_l))\rightarrow \bar P,\bar{p}_0)$ can factor through $V$: $\widetilde{\mu}^0_{0}=(\pi_1(\bar M\times V^0,(\bar{m},v_l))\rightarrow\pi_1(\bar M\times V,(\bar{m},v_l))\xrightarrow{\widetilde{\mu}_{0}} \bar P)$ for some $\widetilde{\mu}_{0}$, lift $(\widetilde{\mu}_{0},\bar p_0)$ first to get $(\widetilde{\psi}_0,p_0)$ with restriction $\widetilde{\psi}_0^0$ on $\bar M\times V^0$.

Paragraph (\ref{thm2}.3) in the proof of Theorem~\ref{thm2} carries over with some obvious modification involving the property of factoring through $V$.
\end{proof}

Below is the last of the 3 steps in constructing the fine moduli space in Proposition~\ref{ram mdli-3}.
\begin{pro}\label{ram mdli-3}
There is a fine moduli space representing $F^T_{U,G}$, the functor for pointed $G$-covers of $(U,u_g)$ tamely ramified over finitely many closed points $T$ on $U$, which is a disjoint union of finitely many ind affine spaces.
\end{pro}
\begin{proof}
 By the same argument for Theorem~\ref{thm1} with slight modification, $F^T_{U,G}=\amalg_{V_l}F^{\rho_{n_l}\bullet/T}_{V_l^0,P}$ which has a fine moduli space by Lemma~\ref{ram mdli-2}.
\end{proof}

Above is the construction of the global side of the local-global principal Proposition~\ref{local-global moduli}, whose local side is the local parameter space in Proposition~\ref{local para space-3}, constructed below in similar 3 steps.

Below are necessary settings of Proposition~\ref{weak local moduli for H}.

\begin{nta}\label{V0t}
Recall $U_0=Spec(k((x)))$ and point $U_0$ at $u_0$. Let $(V_{0t},v_{0t})$ run over all the finitely many connected pointed $\mathbb{Z}/n_t$-covers of $(U_0,u_{0})$, where $n_t$ can be any factor of $n$. Let $V_{0}$ be a connected $\mathbb{Z}/{n}$-cover of $U_0$ given by $k((x))[Y]/(Y^n-x)=k((y))$ with $\bar1 \in \mathbb{Z}/n$ acting on $k((y))$ as $y\mapsto \zeta_n y$. Since $\mathbb{Z}/{n}$ is abelian, $V_0$ can be pointed at any $v_{0}$ over $u_0$, by Remark~\ref{abel}.
\end{nta}
\begin{dfn}\label{w}
Two pointed $Gr$-covers of $(S\times U_0,(s_0,u_0))$ with $(S,s_0)\in \mathcal{S}_1$ are \textit{w-equivalent} if they become isomorphic pulled back to $(\widetilde{T}^0,\tilde{t}_0)$, which is the restriction over $S\times U_0$ of some finite etale cover $(\widetilde{T},\tilde{t}_0)
\rightarrow(S\times \bar U_0,(s_0,u_0))$.

Let $\tilde{\varphi}_i:\pi_1(S\times U_0,(s_0,u_0))\rightarrow Gr$ ($i$=1,2) be two group homomorphisms. They are \textit{w-equivalent} if their corresponding pointed $Gr$-covers of $(S\times U_0,(s_0,u_0))$ are w-equivalent. Denote the w-equivalence class of $\tilde{\varphi}_1$ by $[\tilde{\varphi}_1]^w$.

Let $F^{w\rho}_{V_{0},P}$ be the functor: $\mathcal{S}_1\rightarrow$(Sets), $(S,s_0)\mapsto$\{w-equivalence classes of $\rho$-liftable $P$-covers of $S\times V_{0}$ pointed over $(s_0,v_{0})$\}.
\end{dfn}

\begin{rem}\label{wequi}
The definition of w-equivalence is taken from the 2nd paragraph in the proof of Proposition 2.1 in [H80], which is the right definition of equivalence in the local case to make the proof work.
\end{rem}

Below is the building block needed in the first of the 3 steps in constructing the local parameter space in Proposition~\ref{local para space-3}.
\begin{pro}\label{weak local moduli for H}
Let $(V_0,v_0)$ be given in Notation~\ref{V0t}. Suppose $(\rho,H)$ is in the case of Theorem~\ref{thm3}. Then there exists a fine moduli space $M^{w\rho}_{V_0,H}$ for $F^{w\rho}_{V_0,H}$, the functor for w-equivalence classes of pointed $\rho$-liftable $H$-covers of $(V_0,v_0)$.
\end{pro}
\begin{proof}
The proof is parallel to that of Theorem~\ref{thm3}.

Similarly to the proof of Theorem~\ref{thm3}, start with a short exact sequence $k((y))\xrightarrow{\wp}k((y))\xrightarrow{\pi}
H^1(V_0,H)\rightarrow0$ given by the Artin-Schreier sequence. It can be simplified to $y^{-1}k[y^{-1}]\xrightarrow{\wp}y^{-1}k[y^{-1}]
\xrightarrow{\pi}H^1(V_0,H)\rightarrow0 $ (\ref{weak local moduli for H}.1).

Denote by $\sigma_0$ the automorphism in $Gal(k((y))/k((x)))$ given by $y\mapsto \zeta_n y$. The action of $\rho(-\bar 1)$ on $H$ is given by multiplication by some $e_{\rho}\in \mathbb{F}_q$ ($q=|H|$; see proof of Theorem~\ref{thm3}). Similarly let $D_0$ be the $x^{-1}k[x^{-1}]$-module endomorphism $\sigma_0-e_\rho$ of $y^{-1}k[y^{-1}]$.

Similarly extract from (\ref{weak local moduli for H}.1) an $\mathbb{F}_q$-vector space short exact sequence $KerD_0\xrightarrow[]{\wp}KerD_0\xrightarrow[]{\pi}\mathbb{X}_0\rightarrow0$, where $\mathbb{X}_0$ is the set of $\rho$-liftable pointed $H$-covers of $(V_0,v_0)$. From the $Ker$ exact sequence construct the fine moduli space $M^{w\rho_0}_{V_0,H}$ same as before, which is also an ind affine space. Choose basis $K_{0i}$, an analogue to $K_i$ in the proof of Theorem~\ref{thm3}, inductively for $i\in\mathbb{N}$. The affine space $Spec(k[K_{0i+1}\spcheck-K_{0i}\spcheck])$ can be identified with the $(i+1)$-th piece of $M^{w\rho}_{V_0,H}$; the transition morphism from $Spec(k[K_{0i}\spcheck-K_{0i-1}\spcheck])$ to $Spec(k[K_{0i+1}\spcheck-K_{0i}\spcheck])$ is given by Frobenius as before. Finally, with slight modification to the last two paragraghs of the proof of Theorem~\ref{thm3}, $M^{w\rho}_{V_0,H}$ can be shown to represent $F^{w\rho}_{V_0,H}$.
\end{proof}

\begin{rem}\label{loc uni fml}
Similar to Remark~\ref{uni fml}, a canonical universal family representative over $M^{w\rho}_{V_0,H}$ can be given by $z^q-z=\sum_{\verb"k"_{0i}\in K_{0n}-K_{0n-1}} \verb"k"_{0i}\spcheck\otimes \verb"k"_{0i}$ ($n\geq1$). Precise description can be got by some obvious replacement of symbols in Remark~\ref{uni fml}.
\end{rem}

\begin{rem}\label{l=g bb}
The remark is the base case in the proof for the local-global principal Proposition~\ref{local-global moduli}. Let $\mathbb{A}^{1'}=Spec(k[x^{-1}])$. Suppose $\mathbb{A}^{1'}$ is pointed at $a_g$ such that the map $U_0\rightarrow\mathbb{A}^{1'}$ sends $u_0$ to $a_g$.

Let $V\rightarrow \mathbb{A}^{1'}$ be the $\mathbb{Z}/n$-cover given by $k[x^{-1}][Y^{-1}]/((Y^{-1})^n-x^{-1})=k[y^{-1}]$ ramified at $\infty$, with $\bar1 \in \mathbb{Z}/n$ acting as $y^{-1}\mapsto \zeta_n^{-1} y^{-1}$. Point $V$ at $v_g$ such that $V\rightarrow\mathbb{A}^{1'}$ sends $v_g$ to $a_g$. Its restriction (pullback) at $0$ gives $(V_0,v'_0)$ and let $v_0$ above be $v'_0$:
 \[
 \xymatrix{
  (V_0,v'_0) \ar[d]_{} \ar[r]^{} & (V,v_g)\ar[d]^{} \\
  (U_0,u_0) \ar[r]^{} & (\mathbb{A}^{1'},a_g) .}
 \]
 The constructions show that $M^{w\rho}_{V_0,H}=M^{\rho}_{V,H}$:\\

  The short exact sequence $y^{-1}k[y^{-1}]\xrightarrow{\wp}y^{-1}k[y^{-1}]\xrightarrow{\pi}
  H^1(V_0,H)\rightarrow0$ in the proof of Proposition~\ref{weak local moduli for H}, is similar to the one $k[y^{-1}]\xrightarrow{\wp}k[y^{-1}]\xrightarrow{\pi}
  H^1(V,H)\rightarrow0$ for $V$ in the proof of Theorem~\ref{thm3}, after modding $k[y^{-1}]$ by $k$. The short exact sequence $KerD_0\xrightarrow[]{\wp}KerD_0\xrightarrow[]{\pi}\mathbb{X}_0\rightarrow0$ above, is similar to the short exact sequence $KerD\xrightarrow[]{\wp}KerD\xrightarrow[]{\pi}\mathbb{X}\rightarrow0$ for $V$ in the proof of Theorem~\ref{thm3} and $KerD_0=KerD$. Then the constructions of the two moduli spaces out of the $Ker$ short exact sequences are the same, which shows that $M^{w\rho}_{V_0,H}$ is the same ind scheme as $M^{\rho}_{V,H}$.

Moreover, there is a triangle compatibility diagram. For any pointed $\rho$-liftable $H$-cover $(W,w_g)$ of $(V,v_g)$ corresponding to some $k$-morphism $Spec(k)\xrightarrow{\mathfrak{c}_g}M^{\rho}_{V,H}$, its restriction $(W_0,w_0)$ over $(V_0,v_0)$ is a pointed $\rho$-liftable $H$-cover of $V_0$:
 \[
 \xymatrix{
  (W_0,w_0) \ar[d]_{} \ar[r]^{} & (W,w_g)\ar[d]^{} \\
  (V_0,v_0) \ar[r]^{} & (V,v_g) .}
 \]
 The local cover corresponds to some $k$-morphism $Spec(k)\xrightarrow{\mathfrak{c}_0}M^{w\rho}_{V_0,H}$. The following diagram commutes:
\[\xymatrix{
  Spec(k) \ar[dr]_{\mathfrak{c}_g} \ar[r]^{\mathfrak{c}_0}
                & M^{w\rho}_{V_0,H} \ar[d]^{=}  \\
                & M^{\rho}_{V,H}         .}\eqno{(\ref{l=g bb}.1)} \]

\end{rem}

\begin{rem}\label{alltamecovs}
All the $\mathbb{Z}/n_t$-covers of $U_0$, where $n_t$ can be any factor of $n$, correspond bijectively to all the $\mathbb{Z}/n_t$-covers of $\mathbb{A}^{1'}$ ramified at $\infty$, since these covers can be given by explicit equations of the type in Proposition~\ref{weak local moduli for H} and that in Remark~\ref{l=g bb}.
\end{rem}

Below is the first of the 3 steps in constructing the local parameter space in Proposition~\ref{local para space-3} using the building block in Proposition~\ref{weak local moduli for H}.
\begin{dfn}\label{lftloc}
For the local case, a \textit{pointed $\rho$-liftable $P$-cover} of $(V_0,v_0)$ is defined in the obvious similar way to the global case defined in diagram (\ref{glb mdli}.1). Similarly for a \textit{$\rho$-liftable pair}.

The $k$-points of an ind scheme $M$ \textit{parameterize} certain covers, if there is a bijection $\chi$ together given with $M$ between the set of $k$-points of $M$ and the set of these certain covers.
\end{dfn}

\begin{lem}\label{local para space-1}
Suppose $(\rho,H)$ is in the case of Theorem~\ref{thm3}. There exists a parameter space $M^{p\rho,\bullet}_{V_0,H}$, a disjoint union of finitely many ind affine spaces, whose $k$-points parameterize (see Definition~\ref{lftloc}) all the $\rho$-liftable pairs of $(V_0,v_0)$.
\end{lem}
\begin{proof}
Let $S=Spec(k)$ pointed at $s_0$ that is determined by $v_0$, using diagram (\ref{secn&t}.1). Since $M^{w\rho}_{V_0,H}$ represents $F^{w\rho}_{V_0,H}$, there is a bijection $\chi^{w\rho}_{V_0,H}$ between $F^{w\rho}_{V_0,H}(Spec(k),s_0)$ and $M^{w\rho}_{V_0,H}(Spec(k))$. $F^{w\rho}_{V_0,H}(Spec(k),s_0)$ is the set of pointed $\rho$-liftable $H$-covers of $(V_0,v_0)$.

Let $M^{p\rho,\bullet}_{V_0,H}=\amalg_h M^{w\rho }_{V_0,H}$, an analogue of Corollary~\ref{disj uni over H}. Depending on $(\rho ,H)$, $h$ runs over $H$ or it is just $1$. By the same kind of argument of Corollary~\ref{disj uni over H}, there is a bijection $\chi^{p\rho,\bullet}_{V_0,H}$ between the set of $\rho$-liftable pairs of pointed $H$-covers of $(V_0,v_0)$, and $M^{p\rho ,\bullet}_{V_0 ,H}(Spec(k))$.
\end{proof}

Here is the 2nd of the 3 steps in constructing the local parameter space in Proposition~\ref{local para space-3}.
\begin{lem}\label{local para space-2}
There exists a parameter space $M^{p\rho,\bullet}_{V_0,P}$, a disjoint union of finitely many ind affine spaces, whose $k$-points parameterize (see Definition~\ref{lftloc}) all the $\rho$-liftable pairs of $(V_0,v_0)$.
\end{lem}
\begin{proof}
The proof is parallel to that of Theorem~\ref{thm2} but simpler. It simply replaces some symbols in and do a little modification to the proof of Theorem~\ref{thm2}.

First of all there is no longer an $F$, instead there is $C^{\rho\bullet}_{V_0,P}$ the set of $\rho$-liftable pairs (of pointed $P$-covers) of $(V_0,v_0)$.

In paragraph (\ref{thm2}.1) replace $M^{\bar \rho\bullet}_{V,\bar P}$ by $M^{p\bar\rho\bullet}_{V_0,\bar P}$ and $M^{\rho_0\bullet}_{V,H}$ by $M^{p\rho_0\bullet}_{V_0,H}$.

In paragraph (\ref{thm2}.2) replace $V$ by $V_0$.

In paragraph (\ref{thm2}.3), there is no longer an $S$. Replace every $V$ by $V_0$, and $v_g$ by $v_0$. Replace $(\widetilde \phi,p_1)$ by an element
$(\varphi:\pi_1(V_0,v_0)\rightarrow P,p_1)$ in $C^{\rho\bullet}_{V_0,P}$ and $(\bar{\widetilde \phi},\bar p_1)$ by $(\bar\varphi:\pi_1(V_0,v_0)\rightarrow \bar P,\bar p_1)$. Then replace $\beta$ by a $k$-morphism $Spec(k)\xrightarrow{\mathfrak{c}_\beta} \bar{M}$ and ${\tilde\beta}_*$ by $\tilde{\mathfrak{c}}_{\beta*}$. There is no need for etale descent now and one directly gets a $\mathfrak{c}_\alpha:Spec(k)\xrightarrow{} M^0$. Then replace $M_{V,P}^{\rho\bullet}$ by $M_{V_0,P}^{p\rho\bullet}$. Finally the assignment $\varphi\mapsto(\mathfrak{c}_\alpha,\mathfrak{c}_\beta)$ is a bijection between $C^{\rho\bullet}_{V_0,P}$ and $M(Spec(k))$, which gives the bijection $\chi^{p\rho,\bullet}_{V_0,P}$ desired.
\end{proof}
 Here is the last of the 3 steps in constructing the local parameter space in Proposition~\ref{local para space-3}.
\begin{pro}\label{local para space-3}
There exists a parameter space $M^p_{U_0,G}$, a disjoint union of finitely many ind affine spaces, whose $k$-points parameterize (see Definition~\ref{lftloc}) all the pointed $G$-covers of $(U_0,u_0)$.
\end{pro}
\begin{proof}
Let $M^p_{U_0,G}=\amalg_{V_{0t}}
M^{p\rho_{n_t},\bullet}_{V_{0t},P}$, an analogue of Theorem~\ref{thm1}. Using the argument of Theorem~\ref{thm1} with some obvious modification and Lemma~\ref{local para space-2}, there is a bijection $\chi^p_{U_0,G}$ between $k$-points of $M^p_{U_0,G}$ and pointed $G$-covers of $(U_0,u_0)$.
\end{proof}

Let $M^{\infty}_{\mathbb{A}^{1'},G}$, with $\mathbb{A}^{1'}$ defined in Remark~\ref{l=g bb}, be the short hand notation for $M^{\{\infty\}}_{\mathbb{A}^{1'},G}$, the fine moduli space for $F^{\{\infty\}}_{\mathbb{A}^{1'},G}$ given by Proposition~\ref{ram mdli-3}. Below is the local-global principal that involves the global moduli space in Proposition~\ref{ram mdli-3} and the local parameter space in Proposition~\ref{local para space-3}.
\begin{pro}\label{local-global moduli}
The fine moduli space $M^{\infty}_{\mathbb{A}^{1'},G}$, is the same ind scheme as the parameter space $M^p_{U_0,G}$, compatibly with the inclusion of $U_0$ in $\mathbb{A}^{1'}$ (see diagram (\ref{l=g bb}.1)).
\end{pro}

\begin{proof}
In the construction of both spaces, there are similar 3 steps to the global case, i.e. Theorem~\ref{thm3}$\Rightarrow$Theorem~\ref{thm2}$\Rightarrow$
Theorem~\ref{thm1}. Hence the equality wanted will be proven in similar 3 steps. The bijections $\chi$'s given in Lemma~\ref{local para space-1}, Lemma~\ref{local para space-2} and Proposition~\ref{local para space-3} will be used but not written out unnecessarily.

First both spaces have as building blocks an analogue of the moduli space in Theorem~\ref{thm3}. Let $V_0$ and $V$ be the same as in Remark~\ref{l=g bb}, which shows that the local building block is the same as the global one and a triangle compatibility diagram holds. Then using the same kind of argument as in Corollary~\ref{disj uni over H} $M^{p\rho\bullet}_{V_0,H}=M^{\rho\bullet}_{V,H}$ and a triangle compatibility diagram similar to that in Remark~\ref{l=g bb} holds. Moreover, by Remark~\ref{l=g bb}, Remark~\ref{uni fml}, Corollary~\ref{disj uni over H} and Remark~\ref{loc uni fml}, the canonical $\rho$-liftable universal family representative of $H$-covers of $V_0$ over each connected component of $M^{p\rho\bullet}_{V_0,H}$, is the restriction of the canonical $\rho$-liftable universal family representative of $H$-covers of $V$ over the corresponding connected component of $M^{\rho\bullet}_{V,H}$, which is $M^{\rho\bullet/\{\infty\}}_{V^0,H}$ by Lemma~\ref{ram mdli-1}. Here $U=\mathbb{A}^{1'}$, $T=\{\infty\}$, and $V^0$ is defined at the beginning of this section.

Next $M^{p\rho\bullet}_{V_0,P}=M^{p\rho_0\bullet}_{V_0,H}\times M^{p\bar\rho\bullet}_{V_0,\bar P}$ and $M^{\rho\bullet/\{\infty\}}_{V,P}
=M^{\rho_0\bullet/\{\infty\}}_{V,H}
\times M^{\bar\rho\bullet/\{\infty\}}_{V,\bar P}$ given respectively in Lemma~\ref{local para space-2} and Lemma~\ref{ram mdli-2}. By inductive hypothesis $M^{p\bar\rho\bullet}_{V_0,\bar P}=M^{p\bar\rho\bullet/\{\infty\}}_
{V,\bar P}$, a triangle compatibility diagram similar to that in Remark~\ref{l=g bb} holds, and a $\bar\rho$-liftable universal family representative of pointed $\bar P$-covers of $(V_0,v_0)$ over each connected component of $M^{p\bar\rho\bullet}_{V_0,\bar P}$, is the restriction of a $\bar\rho$-liftable universal family representative of pointed $\bar P$-covers of $(V,v_g)$ over the corresponding connected component of $M^{\bar\rho\bullet/\{\infty\}}_{V,\bar P}$. (Strictly speaking, Definition~\ref{unifmlpair} needs to be used and pairs should be dealt with, which however will make the proof unnecessarily longer.) A $ \rho$-liftable lift of the previous representative can be got from the restriction of a $ \rho$-liftable lift of the latter representative. By this fact and the paragraph above $M^{p\rho\bullet}_{V_0,P}=
M^{\rho\bullet/\{\infty\}}_{V,P}$ and a triangle compatibility diagram similar to that in Remark~\ref{l=g bb} holds.

Finally by Remark~\ref{alltamecovs}, Proposition~\ref{ram mdli-3} and Proposition~\ref{local para space-3}, the proposition follows and there is a triangle compatibility diagram similar to that in Remark~\ref{l=g bb}.
\end{proof}

\begin{cor}\label{local-global}
Any pointed $G$-cover of $(U_0,u_0)$ extends uniquely to a pointed $G$-cover of $(\mathbb{A}^{1'},a_g)$ which is tamely ramified at $\infty$.
\end{cor}
\begin{proof}
By the compatibility assertion in Proposition~\ref{local-global moduli}.
\end{proof}

Here are some necessary settings for the last result in Section~\ref{l vs g}, Proposition~\ref{fin etl-1}.
\begin{nta}\label{chemin}
Let $U_{0_i}$ be the spectrum of the fraction field of the complete local ring at the $i$-th closed point of $\bar U-U$, which is an infinitesimal neighborhood of that point. Let $n'$ be a factor of $n$. Let $V_{0_in'}$ be a fixed connected $\mathbb{Z}/n'$-cover of $U_{0_i}$. All the connected $\mathbb{Z}/n'$-covers of $U_{0_i}$ are isomorphic to $V_{0_in'}$ (see Remark~\ref{alltamecovs}); they only differ by the action of $\mathbb{Z}/n'$. Any two actions differ by an element in $Aut(\mathbb{Z}/n')$.

For any $(n_i)_i$, where $n_i$ is a factor of $n$, there exists a possibly ramified connected $\mathbb{Z}/n$-cover $V$ of $U$ that may ramify at a finite set of closed points $T$ on $U$, such that its ramification index at $U_{0_i}$ is $n_i$. The cover $V$ can be obtained as follows: Suppose $U=Spec(A)$ and denote the fraction field of $A$ by $K$. Pick $a_0\in A$, such that $a_0$ has poles $\sum_i N_iQ_i$, where $Q_i$ is the $i$-th closed point of $\bar U-U$ and $N_i>>0$ with $(N_i,n)=n/n_i$. By Riemann-Roch, such an $a_0$ exists. By adding a constant in $k$ to $a_0$, $Y^n-a_0$ can be assumed an irreducible polynomial in $K[Y]$. Denote $K[Y]/(Y^n-a_0)=K(y)$ by $F$. The normalization of $U$ in $F$ gives $V$, which may ramify over the zeros of $a_0$ on $U$.

Suppose $U_{0_i}$ is pointed at $u_{0_i}$ and $(U_{0_i},u_{0_i})$ maps to $(U,u_{gi})$. Choose $v_{gi}$ such that $(V,v_{gi})\rightarrow(U,u_{gi})$. Let the pointed connected component of $V$'s restriction (pullback) over $U_{0_i}$ be $(V_{i0},v_{i0})$:
 \[
 \xymatrix{
  (V_{i0},v_{i0}) \ar[d]_{} \ar[r]^{} & (V,v_{gi})\ar[d]^{} \\
  (U_{0_i},u_{0_i}) \ar[r]^{} & (U,u_{gi}) .}
 \]
 Then $V_{i0}$ is isomorphic to $V_{0_in_i}$ that is one of those fixed above, as covers of $U_{0_i}$. Choose $v_g$ such that $(V,v_g)\rightarrow(U,u_g)$ and chemins $\omega_i$ from $v_{gi}$ to $v_g$ that induce chemins $\varpi_i$ from $u_{gi}$ to $u_g$.

Let $U^0=U-T$ and $V^0$ be $V$'s restriction over $U^0$, same as the beginning of this section. A $\rho$-liftable pointed $P$-cover of $(V^0,v_g)$ gives a  $\rho_{n_i}$-liftable (see proof of Theorem~\ref{thm1}) pointed $P$-cover of $(V_{i0},v_{i0})$ for each $i$ using the following diagram:
\[\xymatrix{
  \pi_1(V_{i0},v_{i0}) \ar@{^{(}->}[d]_{} \ar[r]^{} &\pi_1(V^0,v_{gi})\ar@{^{(}->}[d]_{} \ar[r]^{\tau_{\omega_i}} & \pi_1(V^0,v_g)\ar@{^{(}->}[d]_{} \ar[r]^{} & P \ar@{^{(}->}[d]^{} \\
  \pi_1(U_{0_i},u_{0_i})\ar[r]^{}&\pi_1(U^0,u_{gi}) \ar[r]^{\tau_{\varpi_i}} &\pi_1(U^0,u_g) \ar[r]^{} & P\rtimes_{\rho}\mathbb{Z}/n ,}
\]
where $\tau_{\omega_i}$ is the isomorphism induced by the chemin $\omega_i$ and similarly for $\tau_{\varpi_i}$.
\end{nta}

Here is a definition involved in the statement of Proposition~\ref{fin etl-1}.
\begin{dfn}\label{def essen surj}
For every $(V_{0n_i},v_{0i})$ a degree $n_i$ cover of $(U_{0_i},u_{0_i})$, denote by $M^{p\rho_{n_i}}_{V_{0n_i},P}$ a connected component (see Remark~\ref{compo}) of $M^{p\rho_{n_i}\bullet}_{V_{0n_i},P}$.

Let $i$ be the index for the $i$-th closed point of $\bar U-U$ and $(n_i)_i$ the same notation in Notation~\ref{chemin}. A morphism from an ind scheme that is a disjoint union of finitely many ind affine spaces, to $\Pi_i M^p_{U_{0_i},G}$ is \textit{essentially surjective}, if for any $(n_i)_i$ there is a connected component of the source ind scheme that maps surjectively (see Definition~\ref{indschmorphdfn} d) to a connected component of the target ind scheme, whose $i$-th factor for each $i$ is $M^{p\rho_{n_i}}_{V_{0n_i},P}$ for some $(V_{0n_i},v_{0i})$ a degree $n_i$ cover of $(U_{0_i},u_{0_i})$.
\end{dfn}

\begin{rem}
The definition of essentially surjective is needed because: Suppose $(V_{i0})_i$ is a tuple whose $i$-th component is the restriction of a possibly ramified $\mathbb{Z}/n$-cover $V$ of $U$ and of degree $n_i$ over $U_{0_i}$. The Galois actions of $\mathbb{Z}/n_i$'s on the $V_{i0}$'s are related to each other as shown in the example below. Thus not every tuple $(V_{0n_i})_i$ (same notation as in Definition~\ref{def essen surj}) could be the image of the restrictions of some $V$. Hence the restriction morphism in Proposition~\ref{fin etl-1} below is not surjective. However in some sense it is surjective, which motivates the definition of essential surjectivity.

Suppose $\textsf{p}=3$. The $\mathbb{Z}/3$-cover of $U=Spec(k[x,x^{-1}])$, the affine line with 0 deleted, given by $V=Spec(k[x,x^{-1}][Y]/(Y^3-x))=Spec(k[x,x^{-1}][y])$ with $\bar 1 \in \mathbb{Z}/3$ acting on $V$ over $U$ as $y\mapsto \zeta_3y$, has restrictions at 0 and $\infty$. At 0, its restriction is a $\mathbb{Z}/3$-cover of $Spec(k((x)))$ given by $V_{00}=Spec(k((x))[Y]/(Y^3-x))=Spec(k((x))[y])$ with $\bar 1 \in \mathbb{Z}/3$ acting on $V_{00}$ over $Spec(k((x)))$ as $y\mapsto \zeta_3y$. At $\infty$, its restriction is a $\mathbb{Z}/3$-cover of $Spec(k((x^{-1})))$ given by $V_{0\infty}=Spec(k((x^{-1}))[Y]/(Y^3-x))=Spec(k((x^{-1}))[y])$ with $\bar 1 \in \mathbb{Z}/3$ acting on $V_{0\infty}$ over $Spec(k((x^{-1})))$ as $y\mapsto \zeta_3y$.

Changing the $\mathbb{Z}/3$ actions on the two local $\mathbb{Z}/3$-covers at 0 and $\infty$ above, the pair of local $\mathbb{Z}/3$-covers (~~($Spec(k((x))[Y]/(Y^3-x))\rightarrow Spec(k((x))),\bar 1:y\mapsto \zeta_3y$), ($Spec(k((x^{-1}))[Y]/(Y^3-x))\rightarrow Spec(k((x^{-1}))),\bar 1:y\mapsto \zeta_3^{-1}y$)~~) got can not come from restrictions of a global $\mathbb{Z}/3$-cover of $Spec(k[x,x^{-1}])$.

\end{rem}

Below is another ingredient involved in the statement of Proposition~\ref{fin etl-1}.

For any $(n_i)_i$, as shown in Notation~\ref{chemin}, there exists a $V_{(n_i)_i}$ that may ramify at a finite set of closed points on $U$, denoted by $T_{V_{(n_i)_i}}$, such that its ramification index at $U_{0_i}$ is $n_i$. Let $T=\cup_{(n_i)_i} T_{V_{(n_i)_i}}$.\\

The last ingredient involved in the statement of Proposition~\ref{fin etl-1}, the restriction morphism, is given in two steps in Lemma~\ref{res mor-1} and Lemma~\ref{res mor-2}. First a map $\texttt{r}$ involved in the statements of Lemma~\ref{res mor-1} and Lemma~\ref{res mor-2}, is defined.

A pointed $G$-cover of $(U^0,u_g)$ gives a local cover of $(U_{0_i},u_{0_i})$ for each $i$: $\pi_1(U_{0_i},u_{0_i})\rightarrow\pi_1(U^0,u_{gi}) \xrightarrow{\tau_{\varpi_i}} \pi_1(U^0,u_g) \rightarrow G$.
Thus there is a map $\texttt{r}$ from the closed points (same as $k$-points) of $M^T_{U,G}$ (see Proposition~\ref{ram mdli-3}), which parameterize certain pointed $G$-covers of $(U^0,u_g)$, to the closed points of $\Pi_i M^p_{U_{0_i},G}$, which parameterize tuples each of which consists of covers indexed by $i$ with the $i$-th entry a pointed $G$-cover of $(U_{0_i},u_{0_i})$. Similarly there is a map $\texttt{r}_0$ from the closed points of $M^{\rho}_{V,H}$ to those of $\Pi_i M^{p\rho_{n_i}}_{V_{i0},H}$.

\begin{lem}\label{res mor-1}
Suppose $(\rho,H)$ is in the case of Theorem~\ref{thm3}. With the notations above, there is a restriction morphism $r_0:M^{\rho}_{V,H}\rightarrow\Pi_i M^{p\rho_{n_i}}_{V_{i0},H}$ such that every closed point of $M^{\rho}_{V,H}$ maps to the same closed point of $\Pi_i M^{p\rho_{n_i}}_{V_{i0},H}$ under $r_0$ or $\texttt{r}_0$.
\end{lem}
\begin{proof}
$r_0$ is given by giving for every $i$ its $i$-th factor using $M^{p\rho_{n_i}}_{V_{i0},H}=M^{w\rho_{n_i}}_{V_{i0},H}$ is a fine moduli space.

Denote by $\widetilde{Z}$ the canonical $\rho$-liftable universal family representative of $H$-covers of $V$ over $M^{\rho}_{V,H}$, which corresponds to, for every point $m$ on $M^{\rho}_{V,H}$, some group homomorphism $\pi_1(M^{\rho}_{V,H}\times V,(m,v_g))\rightarrow H$. The composition $\pi_1(M^{\rho}_{V,H}\times V_{i0},(m,v_{i0}))\rightarrow\pi_1(M^{\rho}_{V,H}\times V,(m,v_{gi}))\xrightarrow{\tilde{\tau}_{\omega_i}}\pi_1(M^{\rho}_{V,H}\times V,(m,v_g))\rightarrow H$, where $\tilde{\tau}_{\omega_i}$ is induced by $\omega_i$ similar to $\tau_{\omega_i}$ given at the end of Notation~\ref{chemin}, gives a pointed $H$-cover of $(M^{\rho}_{V,H}\times V_{i0},(m,v_{i0}))$ the non pointed version of which is denoted by $\widetilde{Z}_0$. By Remark~\ref{abel}, base points here do not matter. It is the restriction (pullback) of $\widetilde{Z}$ to $M^{\rho}_{V,H}\times V_{i0}$:
\[\xymatrix{
  \widetilde{Z}_0 \ar[d]_{} \ar[r]^{} & \widetilde{Z}\ar[d]^{} \\
  M^{\rho}_{V,H}\times V_{i0} \ar[r]^{} &   M^{\rho}_{V,H}\times V .}
\]
Since $M^{p\rho_{n_i}}_{V_{i0},H}=M^{w\rho_{n_i}}_{V_{i0},H}$ and $M^{w\rho_{n_i}}_{V_{i0},H}$ represents $F^{w\rho_{n_i}}_{V_{i0},H}$ by Proposition~\ref{weak local moduli for H}, there is a morphism $M^{\rho}_{V,H}\xrightarrow{r_{0i}} M^{p\rho_{n_i}}_{V_{i0},H}$ given by $\widetilde{Z}_0$. A different base point $m'$ gives the same $r_{0i}$. Then define $r_0:=(r_{0i})_i$.

Now it is enough to verify that a closed point $m'$ of $M^{\rho}_{V,H}$ maps to the same closed point under $r_{0i}$ or $\texttt{r}_{0i}$, where $\texttt{r}_{0i}$ is the $i$-th factor of $\texttt{r}_0$. Tracking definitions, $r_{0i}(m')$ represents the restriction to $(V_{i0},v_{i0})$ of the pointed $H$-cover of $(V,v_g)$ represented by $m'$. And $\texttt{r}_i$ does the same thing by its definition. So $r_{0i}$ and $\texttt{r}_{0i}$ agree.
\end{proof}

\begin{lem}\label{res mor-2}
Let $G=H\rtimes_{\rho}\mathbb{Z}/n$ for some $(\rho,H)$ in the case of Theorem~\ref{thm3}. There is a restriction morphism $r: M^T_{U,G}\xrightarrow{}\Pi_i M^p_{U_{0_i},G}$, such that every closed point of $M^T_{U,G}$ maps to the same closed point of $\Pi_i M^p_{U_{0_i},G}$ under $r$ or $\texttt{r}$, where $\texttt{r}$ is defined above Lemma~\ref{res mor-1}.
\end{lem}
\begin{proof}
By construction, $M^T_{U,G}$ and $\Pi_i M^p_{U_{0_i},G}$ are both a disjoint union of finitely many ind affine spaces. The morphism $r$ will be given for each connected component of $M^T_{U,G}$.

Proposition~\ref{ram mdli-3} and Lemma~\ref{ram mdli-1} give $M^T_{U,G}=\amalg_{V_l}M^{\rho_{n_l}\bullet/T}_{V_l^0,H}$ and $M^{\rho_{n_l}\bullet/T}_{V_l^0,H}=\amalg_h M^{\rho_{n_l}}_{V_l,H}$ respectively. A connected component of $M^T_{U,G}$ is of the form $M^{\rho_{n_l}}_{V_l,H}$.

Let the pointed connected component of $V_l$ over $U_{0_i}$ be $V_{li0}$, a $\mathbb{Z}/n_i$-cover of $U_{0_i}$. Using Notation~\ref{chemin}, $r$ should map $M^{\rho_{n_l}\bullet/T}_{V_l^0,H}$  to $\Pi_i M^{p\rho_{n_i}\bullet}_{V_{li0},H}$, since it is required to agree with the map $\texttt{r}$ on closed points. Similarly the target connected component of each connected component of $M^{\rho_{n_l}\bullet/T}_{V_l^0,H}$ under $r$ can be identified. Denote by $ M^{\rho_{n_l}}_{V_l,H }$ a connected component of $M^{\rho_{n_l}\bullet/T}_{V_l^0,H}$, by $\Pi_i M^{p\rho_{n_i}}_{V_{li0},H }$ its target connected component, and by $r_{lj}$ (suppose $ M^{\rho_{n_l}}_{V_l,H }$ is the $j$-th component of $M^{\rho_{n_l}\bullet/T}_{V_l^0,H}$) the restriction of $r$ on $ M^{\rho_{n_l}}_{V_l,H }$.

Finally let $r_{lj}$ be the restriction morphism given in Lemma~\ref{res mor-1} for every index $lj$. One can check that the morphism $r$ satisfies the requirement.
\end{proof}

With the preparation from Notation~\ref{chemin} to Lemma~\ref{res mor-2}, the last result in Section~\ref{l vs g} can be given.
\begin{pro}\label{fin etl-1}
Let $G=H\rtimes_{\rho}\mathbb{Z}/n$ for some $(\rho,H)$ in the case of Theorem~\ref{thm3}. The restriction morphism $M^T_{U,G}\xrightarrow{r}\Pi_i M^p_{U_{0_i},G}$ given in Lemma~\ref{res mor-2} is essentially surjective and finite. And the degrees of $r$ on different connected components of $M^T_{U,G}$ are all powers of $\mathsf{p}$.
\end{pro}
\begin{proof}
The proof follows the points of the proof of Proposition 2.7 in [H80]. A calculation of the dimensions of the $n$-th pieces of both source and target shows that they are the same. By this fact the map $r$ restricted on each connected component of the source can be proven surjective. Then all the three statements follow.

With the same notations as in the proof of Lemma~\ref{res mor-2}, denote a connected component of $M^{\rho_{n_l}\bullet/T}_{V_l^0,H}$ by $ M^{\rho_{n_l}}_{V_l,H}$, whose $n$-th piece is $M^{\rho_{n_l}}_{V_l,H,n}$. Denote the connected component of $\Pi_iM^{p\rho_{n_i}\bullet}_{V_{li0},H}$, which $M^{\rho_{n_l}}_{V_l,H}$ maps to under $r_{lj}$, by $\Pi_iM^{p\rho_{n_i}}_{V_{li0},H}$, whose $n$-th piece is $\Pi_iM^{p\rho_{n_i}}_{V_{li0},H,n}$.

For $n>>0$, Riemann-Roch shows that the dimension of $M^{\rho_{n_l}}_{V_l,H,n}$ is at least that of $\Pi_i M^{p\rho_{n_i}}_{V_{li0},H,n}$ for a subsequence $\{n_\mathsf{k}\}$ of $\mathbb{N}$: A simpler but similar computation is done in the 1st paragraph of the proof of Proposition 2.7 in [H80]. Here pass from $V_l$ to $U$ first using ramification indices and then do a similar computation to [H80]. Denote by $\Sigma_i d_{i0n}$ the dimension of $\Pi_i M^{p\rho_{n_i}}_{V_{li0},H,n}$. For $n>>0$ Riemann-Roch gives $\Sigma_i d_{i0n}=\Sigma_i (\lfloor\frac{q^n-i_0}{n_i}\rfloor-\lfloor\frac{q^{n-1}-i_0}{n_i}\rfloor)$ with some natural number $i_0 $ between 0 and $n_l$. Denoted by $d_n$ the dimension of $M^{\rho_{n_l}}_{V_l,H,n}$. For $n>>0$ a similar computation gives $d_n=\Sigma_i(\lfloor\frac{q^n+\delta_i}{n_i}\rfloor-
\lfloor\frac{q^{n-1}+\delta_i}{n_i}\rfloor)$ for some $\delta_i\in\mathbb{Q}$. Since the remainder of $q^n$ divided by $n_i$ is periodic for $n\in\mathbb{N}$ there is a subsequence  $\{n_\mathsf{k}\}$ of $\mathbb{N}$ such that $d_{n_\mathsf{k}}\geq \Sigma_i d_{i0n_\mathsf{k}}$.

$r_{lj}$ is quasi finite of degree a $\textsf{p}$-power: The restriction of $r_{lj}$ on the $n$-th piece of $M^{\rho_{n_l}}_{V_l,H}$ gives $ M^{\rho_{n_l}}_{V_l,H,n}\xrightarrow{r_{ljn}}\Pi_i
M^{p\rho_{n_i}}_{V_{li0},H,n}$ using Lemma~\ref{res mor-1}, which is in fact a homomorphism between $\mathbb{F}_\mathsf{p}$-vector spaces (The closed points of any moduli space involved here form a $\mathbb{F}_\mathsf{p}$-vector space, by the definitions of the moduli spaces.). Since there are, up to isomorphism, only finitely many pointed (etale) $P$-covers of the completion $X_l=\bar V_l$ [SGA1, X 2.12], the kernel of $r_{ljn}$ is finite (and equal to this number when $n>>0$). Thus every non empty fiber of $r_{ljn}$ (hence of $r_{lj}$) contains the same finite number of points. This number is a power of $\mathsf{p}$, being the cardinality of a $\mathbb{F}_\mathsf{p}$-vector space.

The 2nd paragraph in the proof of Proposition 2.7 in [H80] shows that since $d_{n_\mathsf{k}}\geq \Sigma_i d_{i0n_\mathsf{k}}$ for every $\mathsf{k}$ large enough and $r_{ljn_\mathsf{k}}$ is quasi finite, the morphism $M^{\rho_{n_l}}_{V_l,H,n_\mathsf{k}}\rightarrow\Pi_i M^{p\rho_{n_i}}_{V_{li0},H,n_\mathsf{k}}$ is surjective. Thus for every $n$ large enough, using Lemma~\ref{res mor-1}, the morphism $M^{\rho_{n_l}}_{V_l,H,n}\rightarrow\Pi_i M^{p\rho_{n_i}}_{V_{li0},H,n}$ is surjective. Hence the map $r$ restricted on every connected component of $M^{\rho_{n_l}\bullet/T}_{V_l^0,H}$ maps surjectively to a connected component of $\Pi_i M^{p\rho_{n_i}\bullet}_{V_{li0},H}$.

Direct computation shows that the restriction morphism is finite. The choice of $T$ shows that $r$ is essentially surjective.
\end{proof}

\textbf{References}

[H80] D. Harbater. Moduli of $\mathsf{p}$-covers of curves. Commu. Algebra, 8(12): 1095-1122, 1980

[K86] N. Katz. Local-to-global extensions of representations of fundamental groups. Ann. Inst. Fourier (Grenoble), 36(4):69--106, 1986.

[M13] J. Milne. Lectures on Etale Cohomology v2.21. On www.jmilne.org. 2013

[P02] R. Pries. Families of wildly ramified covers of curves. Amer. J. Math. 124 (2002), 4:737-768.

[Serre] J. -P. Serre. Cohomologie Galoisienne. Springer, New York, 1964

[SGA 1] A. Grothendieck, et. al. Revetements etales et groupe fondamental. LNM 224, Springer, New York, 1971.

[S09] T. Szamuely. Galois groups and fundamental groups. Cambridge university press, 2009.

[TY17] F. Tonini and T. Yasuda. Moduli of formal torsors. arXiv:1709.01705v1\\

\noindent {\bf Author Information:}

jmor.zhang47@foxmail.com
\medskip

\end{document}